\documentclass[10pt]{article}
\usepackage[hmargin=3cm,vmargin=3cm]{geometry}
\usepackage{amsmath,amssymb,paralist,url,hyperref,doi}
\usepackage{palatino,epsfig,latexsym,xfrac} 
\usepackage{setspace}
\usepackage{multirow} 
\usepackage{color}
\usepackage{lineno}
\usepackage{comment}
\usepackage{graphicx}
\usepackage{wrapfig}
\usepackage{balance} 
\def\prob{{\mathbf{Pr}}}
\def\E{\mathop{{\mathbb{E}}}}

\usepackage{algorithmic}
\usepackage[linesnumbered,boxed]{algorithm2e} %
\DontPrintSemicolon
\usepackage{authblk}
\usepackage{tikz}
\usepackage{pgfplots}

\usepackage{enumitem,float}
\title{Correlated Mutations for Integer Programming} 
\author{Ofer M. Shir$^{1,2}$\footnote{Corresponding author: \href{mailto:ofersh@telhai.ac.il}{\texttt{ofersh@telhai.ac.il}}.} and Michael Emmerich$^{3}$}
\affil{ \begin{small} $^{1}$ Tel-Hai College, Upper Galilee, Israel\\$^{2}$ Migal Institute, Qiryat Shemona, Israel\\$^{3}$ University of Jyv{\"a}skyl{\"a}, Finland\end{small}}
\date{}
\begin{document}
\maketitle    
\SetKwInOut{Input}{input}
\SetKwInOut{Output}{output}
\SetFuncSty{texttt}
\SetKwRepeat{Do}{do}{while}
%
\begin{abstract}
Even with the recent theoretical advancements that dramatically reduced the complexity of Integer Programming (IP), heuristics remain the dominant problem-solvers for this difficult category. This study seeks to establish the groundwork for Integer Evolution Strategies (IESs), a class of randomized search heuristics inherently designed for continuous spaces. IESs already excel in treating IP in practice, but accomplish it via discretization and by applying sophisticated patches to their continuous operators, while persistently using the $\ell_2$-norm as their operation pillar. We lay foundations for discrete search, by adopting the $\ell_1$-norm, accounting for the suitable step-size, and questioning alternative measures to quantify correlations over the integer lattice. We focus on mutation distributions for unbounded integer decision variables. We briefly discuss a couple of candidate discrete probabilities induced by the uniform and binomial distributions, which we show to possess less appealing theoretical properties, and then narrow down to the Truncated Normal (TN) and Double Geometric (DG) distributions. We explore their theoretical properties, including entropy functions, and propose a procedure to generate scalable correlated mutation distributions. Our investigations are accompanied by extensive numerical simulations, which consistently support the claim that the DG distribution is better suited for unbounded integer search. We link our theoretical perspective to empirical evidence indicating that an IES with correlated DG mutations outperformed other strategies over non-separable quadratic IP. We conclude that while the replacement of the default TN distribution by the DG is theoretically justified and practically beneficial, the truly crucial change lies in adopting the $\ell_1$-norm over the $\ell_2$-norm.
\end{abstract}
\textbf{Keywords}: Integer Programming; multivariate discrete distributions; Truncated Normal; Double Geometric, evolution strategies; $\ell_1$-norm.

\section{Introduction}
Integer Programming (IP) comprises a rich family of optimization problems with discrete decision variables and is central to both practical computation and the theoretical foundations of Computer Science \cite{PapadimitriouSteiglitz,Schrijver}. Whereas IP is often viewed as a generalization of linear programming for combinatorial optimization, our paper focuses on methods for unconstrained or loosely constrained problems \cite{ShirMOQC}. As opposed to categorical integer variables, we consider metric (or \textit{ordinal}) integer variables -- whose objective‑function values are correlated with the $\ell_{1}$ (generalized Hamming) distance among decision variable vectors. 
Applications include selecting pipeline diameters from a discrete set of options \cite{reehuis2011multiobjective} and determining the number of stages in separation processes within chemical plants \cite{emmerich2000mixed}. Metric integer parameters also arise in hyperparameter tuning for machine learning, image processing, and inventory planning, where decisions must be expressed as integer values \cite{li2013mixed}.

However, tackling such integer optimization problems remains computationally formidable in the worst case.
The \textit{decision} version of IP is $\mathcal{NP}$‑hard, and for decades the fastest exact algorithm was due to
Kannan, with running time \(\text{poly}(n)\,2^{\mathcal{O}(n)}\) for an instance with \(n\) variables~\cite{Kannan1988}.
A recent breakthrough by Rothvoss and Reis tightens this bound: they present a randomized algorithm that solves any IP in \((\log(2n))^{\mathcal{O}(n)}\) steps~\cite{RothvossReis2024FasterIntegerProgramming}.
Despite this major theoretical advancement and the formulation of a faster algorithm, heuristics remain the primary tool to address IP problem solving in practice. In particular, there is a need for solvers that can work in an opaque-box setting, that is, without explicitly knowing the mathematical expression of the objective functions. 
Also in other aspects, unconstrained IP significantly differs from its continuous counterpart.  
For example, a convex quadratic optimization problem can become multi-modal when truncated -- due to discretization effects -- as demonstrated in the Supplementary Material.


When it comes to heuristic methods for solving unconstrained optimization problems, Evolution Strategies (ESs) \cite{Rechenberg71,Schwefel65,Baeck2013contemporary} are effective randomized search heuristics. They stand out from many other evolutionary Algorithms due to search mechanisms motivated by a blend of inspiration from biological evolution and a thorough theoretical analysis \cite{Beyer}, advanced step size adaptation mechanisms \cite{Baeck-book}, and their invariance to monotonic transformations of the objective function as well as rotations in the decision space \cite{InvarianceHansen2000}. 
Although ESs are commonly used in continuous spaces, their paradigms have been adapted to treat integer \cite{Schwefel65,rudolph1994evolutionary} and mixed integer (MI) spaces (see, e.g., \cite{Baeck95MIES,emmerich2000mixed,li2013mixed,moMIQC2024}).  
However, so far, these adaptations primarily rely on independent multivariate distributions or truncated versions of mutation distributions originally designed for continuous settings.
Although integer and MI optimization are areas of significant practical importance, particularly in the context of black-box models, the interest in their theoretical foundations has only recently increased \cite{miBBOB2019,Arnold2023_miConstraints}. Key questions remain unresolved despite modern ESs' capability to solve such problems. 
In fact, modern ESs consistently demonstrate strong capabilities in addressing MI optimization problems (CMA-ES-wM \cite{CMAESwM2022} and CMA-ES-IH \cite{cmaIHgecco2023}) while still operating with the Truncated Normal (TN) distribution under the hood after certain adaptations. 
Having such success in practice has not motivated the question of addressing open theoretical questions and especially investigating the suitability of alternative discrete distributions for integer mutations, among which the effectiveness of the TN distribution remains unclear. 
Rudolph, following the principle of maximum entropy, justified the use of the Double Geometric (DG) distribution for unbounded integer optimization in a singular study published in 1994 \cite{rudolph1994evolutionary}. Rudolph proposed therein a mutation operator for integer search-spaces, grounded in statistical theory, using an independent multivariate sample from the DG distribution. Following it, the DG distribution has been utilized in MIESs \cite{emmerich2000mixed,li2013mixed}, with empirical evidence demonstrating that mutative step-size control can effectively track nearly optimal step-sizes for the so-called Progress Rate \cite{li2013mixed}.
However, such studies with DG-based mutations have remained isolated examples compared to the multivariate TN distribution, which has been performing successfully and, therefore, like the ``winning horse in midstream'', has kept its steady state as the first choice for the designers and practitioners of MIESs. In this context, a recent empirical study \cite{abnormalmutations2025} demonstrated that Gaussianity is not essential for ESs' strong performance over \textit{continuous} spaces, motivating further exploration of alternative distributions.
In addition to this momentum effect, another reason for persisting with the TN distribution might be that a multivariate correlation structure has never been defined for the DG distribution, and it is an open question to which extent such a correlation structure can be imposed. 
This key question will drive us to investigate the foundations of search over the integer lattice, where our working hypothesis states that the $\ell_1$-norm is the pillar of the search. 
Moreover, unlike the continuous domain which enjoys direct $\ell_2$-norm-driven relations between the statistical covariance to translation and rotation transformations, the existence of such relations or even proper definitions of non-trivial rotations over the integer lattice remains an open issue. 
In the context of ESs' research, the required norm-shift from $\ell_2$ to $\ell_1$ constitutes a major challenge with a limited ability to capitalize on the large volume of well-established $\ell_2$-driven results.

\paragraph{Research Questions} 
Motivated by the differences in geometric structure between continuous and integer spaces and by empirical observations on mutation behavior, this paper investigates the following research questions:
\begin{enumerate}
    \item \textbf{Geometry and Entropy of Integer Mutation Distributions.} 
    \begin{itemize}
        \item \textbf{Open Question:} How can mutation operators for integer ESs (IESs) be designed to respect the inherent $\ell_1$-norm (Manhattan) geometry of the integer lattice? How does its entropy and symmetry compare to other proposed distributions, such as the TN distribution?
        \item \textbf{Hypothesis:} For integer optimization problems $\ell_1$-invariant operators are better suited than $\ell_2$-invariant ones. Moreover, the mean $\ell_1$ metric is a more meaningful measure of the average step-size than the standard deviation. 
    \end{itemize}
    \item \textbf{Correlated Integer Mutations via DG Distributions.}
    \begin{itemize}
        \item \textbf{Open Question:} Can a correlated version of the DG distribution be constructed that accurately preserves intended dependencies among integer variables? What is the exact formulation and what are the properties of the TN distribution?
        \item \textbf{Hypothesis:} A sampling scheme based on the DG distribution can be extended to the correlated case while preserving covariance structure and yielding higher entropy than truncated Gaussians. The DG distributions, which can be adjusted to a given covariance structure in a rather straightforward manner, is not well suited to reflect $\ell_1$ -- symmetries, which might make it less suitable for unconstrained integer optimization.
    \end{itemize}
    \item \textbf{Stagnation near Optima: A Discreteness-Induced Phase Transition.}
    \begin{itemize}
        \item \textbf{Open Question:} Why do integer-based evolutionary strategies systematically stagnate near the optimum, and can mutation dynamics be adapted to mitigate this behavior?
        \item \textbf{Hypothesis:} A phase-transition-like stagnation effect occurs due to the reduced move options in the near-optimal integer lattice; standard adaptation strategies fail to overcome it without explicitly accounting for discreteness.
    \end{itemize}
\end{enumerate}

\paragraph{Paper Organization}
We provide the mathematical preliminaries in Section \ref{sec:preliminaries}, where we also state our working hypothesis and present basic numerical evaluations of our mathematical formulations. Section \ref{sec:correlations} is dedicated to the rotation transformation for generating correlated mutations, while Section \ref{sec:entropy} covers the \textit{entropy} perspective of the investigated mutation distributions. 
In Section \ref{sec:experiments} we present practical optimization problem-solving of integer quadratic programs as a proof-of-concept. 
Finally, we summarize our work in Section \ref{sec:summary}.
\section{Preliminaries} \label{sec:preliminaries}
The purpose of this section is to present the necessary preliminaries, including our working hypotheses, and then describe four discrete distributions of interest, namely the Discrete Uniform, the Shifted Binomial, the Truncated Normal and the Double Geometric. 
\subsection{Framework and Notation}
The general framework that underlies our work is MI optimization problem-solving, subscribing to the following formulation: 
\begin{align}\label{eq:generalMI}
\begin{array}{ll}
\displaystyle \textrm{minimize}_{\vec{y}} & f\left( \vec{y} \right)\\
\displaystyle \textrm{subject to:}& g\left(\vec{y}\right)\\
\displaystyle & \vec{y}\in\mathbb{R}^{D}\\ 
\displaystyle & y_i\in \mathbb{Z} \quad \forall i \in I,
\end{array} 
\end{align}
where $f$ and $g$ constitute the objective and constraint functions, respectively. 
Importantly, the $D$-dimensional decision vector $\vec{y}$ is constructed by $n_r$ real-valued decision variables followed by $n_z$ (metric) integer decision variables that are defined by the so-called \textit{index set} $I:=\left\{n_r+1,\ldots,n_r+n_z \right\}$: $\forall i \in I \quad y_i \in \mathbb{Z}$. 
As mentioned earlier, our consideration excludes categorical decision variables, although they may as well exist in practical modeling (see, e.g., \cite{Cucumbers-Postharvest_Shir2022}) -- but exceed the scope herein.\\  
\textbf{While Eq.~\eqref{eq:generalMI} defines the context and induces the applicability and relevance of this study thereafter, we focus entirely on the integer subspace, and thus refer to the integer decision vector as our target, $\displaystyle \vec{x}\in \mathbb{Z}^n$, fully mapped onto the index set $I$ and satisfying $n=n_z$}.

During the solution process of the optimization problem, an existing integer decision vector $\vec{x}_{\text{CURR}}$ undergoes variations to generate a candidate vector $\vec{x}_{\text{NEW}}$:
\[
\displaystyle \vec{x}_{\text{NEW}} = \vec{x}_{\text{CURR}} + \vec{z}.
\]
$\vec{z} \in \mathbb{Z}^n$ is the mutation vector, which is drawn from a certain multivariate distribution. 
This distribution may possess a correlated or uncorrelated nature with respect to the generated random variables (that is, with respect to the independence of $z_i$).\\ 
Unless specified otherwise, we assume here uncorrelated distributions (i.e., we assume $z_i$ to be independent). 

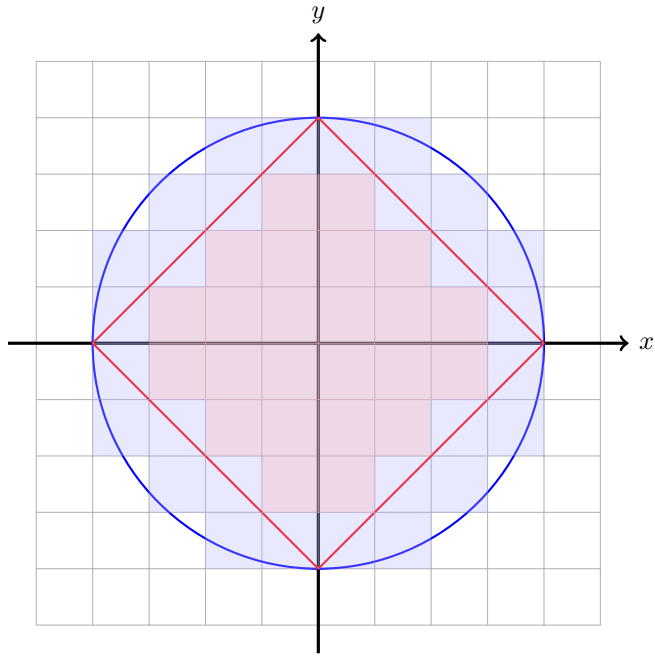
\begin{figure}
    \centering
\begin{tikzpicture}[scale=0.75] 
  \tikzset{
    grid/.style = {step=1cm, thin, gray!60},
    axis/.style = {->, very thick},
    fillE/.style = {blue!30, opacity=0.3},  
    fillL/.style = {red!30,  opacity=0.3}   
  }

  \draw[grid] (-5,-5) grid (5,5);
  \draw[axis] (-5.5,0) -- (5.5,0) node[right] {$x$};
  \draw[axis] (0,-5.5) -- (0,5.5) node[above] {$y$};

  \draw[blue, thick] (0,0) circle[radius=4];

  \draw[red, thick] ( 4,0) -- (0, 4) -- (-4,0) -- (0,-4) -- cycle;

  \foreach \i in {-5,...,4}{
    \foreach \j in {-5,...,4}{
      \pgfmathsetmacro\cx{\i+0.5}
      \pgfmathsetmacro\cy{\j+0.5}
      \pgfmathsetmacro\r{sqrt(\cx*\cx+\cy*\cy)}
      \ifdim \r pt < 4pt
        \fill[fillE] (\i,\j) rectangle ++(1,1);
      \fi
    }
  }

  \foreach \i in {-5,...,4}{
    \foreach \j in {-5,...,4}{
      \pgfmathsetmacro\cx{\i+0.5}
      \pgfmathsetmacro\cy{\j+0.5}
      \pgfmathsetmacro\lone{abs(\cx)+abs(\cy)}
      \ifdim \lone pt < 4pt
        \fill[fillL] (\i,\j) rectangle ++(1,1);
      \fi
    }
  }
\end{tikzpicture}
    \caption{A visual comparison of the $\ell_1$-norm versus the $\ell_2$-norm over the 2D integer lattice using spheres of integer radius 4: an $\ell_1$-sphere (red) versus an $\ell_2$-sphere (blue).}
    \label{fig:ell1}
\end{figure}

To assess the mutation vectors, we utilize the $\ell_1$-norm, which reflects vector length over the integer lattice in the sense that it quantifies the minimal number of elementary steps required to reach one vector from another: 
\begin{align}
    \displaystyle \| \vec{z} \|_1 := \sum_{i=1}^{n} |z_i|.
\end{align}
See Figure \ref{fig:ell1} for an illustrated comparison of the two norms in 2D, where an $\ell_1$-sphere is depicted as opposed to the $\ell_2$-sphere, when visualized over the same integer lattice. The $\ell_1$-norm, in contrast to the $\ell_2$-norm, provides a natural way to describe spheres within the integer lattice. For positive integer radii, the set of points at a fixed $\ell_1$ distance from the origin forms a \textit{discrete sphere} (or ``diamond'') that lies entirely within the integer lattice. These $\ell_1$-spheres are fully representable, meaning that every point on the sphere has integer coordinates, which is generally not true for $\ell_2$-spheres. \\
\textbf{The expected $\ell_1$-norm constitutes the mean step-size}, 
\begin{align}
    \displaystyle S:= \E \left[ \|\vec{z} \|_1 \right] = \sum_{i=1}^{n} \E \left[ |z_i|_1 \right],
\end{align}
due to the stochastic independence. When the $n$ random variables are identically drawn under the same settings, this mean step-size becomes $\displaystyle S = n\cdot \E \left[ |z_1| \right]$.

The probability distributions under consideration are discrete and they are formulated by means of the probability of the mutation element to take an integer value $k$:
\[
\displaystyle \prob \left\{ z = k \right\} = p_k,\quad \sum_k p_k=1.0.
\]

Finally, we mention Shannon's entropy as a measure of averaged information, which quantifies the degree of unpredictability of a random variable's possible outcomes:
\begin{align}\label{eq:Shannon}
  \displaystyle H:= -\sum_{k=-\infty}^\infty p_k \log_2 p_k~.
\end{align}


\subsection{Working Hypothesis}
Our working hypothesis states that \textbf{the $\ell_1$-norm is the natural measure over the integer lattice}, since it directly measures the sum of absolute differences along each coordinate axis. 
At the same time, the $\ell_2$-norm, which considers the Euclidean interpretation, may not be as directly meaningful in discrete lattices, especially in the context of integer mutations.  
We therefore adopt the $\ell_1$-norm and claim that the shift from the $\ell_2$-norm within ESs' existing work necessitates careful consideration of basic components such as the mutative step-size, but it also impacts other considerations, including statistical measures. 
Importantly, the \textit{statistical covariance matrix} is interpretable by means of the $\ell_2$-norm: the variance of a single variable is the square of the $\ell_2$-norm of its deviation from the mean, and the covariance between two variables is the $\ell_2$-norm of their joint deviations. 
As such, the covariance matrix is not necessarily a proper measure to consider in the context of correlations among variables over the integer lattice. 
The $\ell_2$-based Pearson correlation measure also falls into this category, and therefore has a limited potential to statistically describe such variables. 


\subsection{Discrete Mutation Distributions}
Studying the effectiveness of discrete randomized heuristic search in integer search-spaces dates back to the 1960's and 1970's. 
Kelahan and Gaddy proposed to use a bilateral power distribution that adaptively shrank by a geometric schedule \cite{Kelahan1978-RandomSearchMIO}.
Rudolph's early research, dating back to 1994, questioned the suitability of random distributions for evolutionary mutations in unbounded integer spaces \cite{rudolph1994evolutionary}. In particular, Rudolph showed that integer mutations utilizing the geometric distribution possess maximum entropy, and demonstrated the effectiveness of such mutations in practice.
The goal of this section is to present in detail two probability distributions of interest, the TN and DG distributions. 
We begin by mentioning two additional discrete distributions, which are intuitive and may serve as reference distributions. 
\subsubsection{The Discrete Uniform (DU) Distribution}
Playing the role of the most intuitive manner to generate random integers, the discrete uniform distribution of interest draws a random variable from a predefined range at equal probabilities, which in our context of search mutations corresponds to $\left\{-N,~-N+1,\ldots,~0,\ldots,~N-1,~N \right\}$ :
\begin{align}\label{eq:discreteUniform}   
\displaystyle \prob \left\{ X = k \right\} = \frac{1}{2N+1} \quad \forall k \in \left\{-N,-N+1,\ldots,N-1,N \right\}.
\end{align}
This range dictates the variation's strength, and therefore, it is possible to control the mutation step by considering the 
expected $\ell_1$-norm, which reads $\displaystyle S_{DU}= n\cdot \frac{N(N+1)}{2N+1}$. Given a desired mean step-size $S_{*}$, by solving for $N$, it may be set to  
\begin{align}\label{eq:NbyL1normUniform}
    \displaystyle N_{*} := \texttt{int}\left[ \frac{2(\frac{S_{*}}{n}) - 1 + \sqrt{1 + (\frac{S_{*}}{n})^2}}{2} \right].
\end{align}
Finally, the entropy of this distribution is simply 
\begin{align}\label{eq:HDU}
    \displaystyle H_{DU} = \log_2 (2N+1).
\end{align}
\subsubsection{The Shifted Binomial (SB) Distribution}
Describing a certain number of successful trials out of a fixed budget $N$, we consider this discrete probability distribution in our integer context with an equal probability for success/failure ($p=\frac{1}{2}$). Focusing on the symmetric case which enables generation of negative as well as positive values, the random variable is shifted, $Y:=X-\frac{N}{2}$, where \( X \) is the binomial random variable with parameters \( N \) and \( p=\frac{1}{2} \). 
The probability of $Y$ to get a value $k$ is given by 
\begin{align} \label{eq:shiftedBinomial}   
\displaystyle \prob \left\{ Y = k \right\} = \prob \left\{ X = k + \frac{N}{2} \right\} = \binom{N}{k + \frac{N}{2}} p^{k + \frac{N}{2}} (1 - p)^{N - (k + \frac{N}{2})} ~.
\end{align}
Notably, the budget of trials $N$ defines the range of generated numbers, and in our context of integer mutations -- it governs the step's strength.
Equivalently to the DU distribution, it is possible to control the step by considering the expected $\ell_1$-norm. For this distribution, it may be approximated via the so-called Mean Absolute Deviation which yields 
$\displaystyle S_{SB} \approx n\cdot \sqrt{\frac{2N}{\pi}}$. 
Accordingly, given a desired step-size $S_{*}$, the distribution can be tuned accordingly by computing $N_{*}$: 
\begin{align}\label{eq:NbyS-Bimonial}
    \displaystyle N_{*} := \texttt{int}\left[ \frac{\pi}{2} \left(\frac{S_{*}}{n} \right)^2  \right].
\end{align}
However, \textbf{despite this artificial step-size controllablity, the fit of this distribution to unbounded search remains unnatural}. 
\subsubsection{The Truncated Normal (TN) Distribution}
\begin{figure}
    \centering
    \includegraphics[width=0.75\columnwidth]{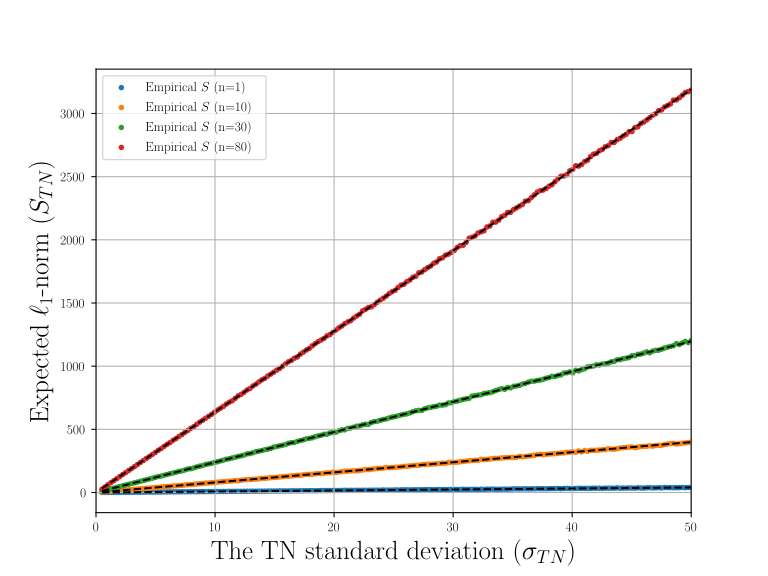}
    \caption{Expected $l_1$-norm $S$ vs. $\sigma$ for the TN Distribution per various vector dimensions: $n\in \{1,10,30,80\}$. The curves depict the empirical expectation (via colored circles) while exercising the approximation of $S_{TN}$ given in Eq. \eqref{eq:TN_S-sigma} (via a dashed line) -- which demonstrate excellent fit.}
    \label{fig:S-vs-sigma_TN}
\end{figure}

The Normal distribution constitutes the continuous limit of the Binomial distribution when the trials' number $n$ increases and becomes sufficiently large.
Given a normally distributed random variable, $z\sim \mathcal{N}(0,\sigma^2)$, which may take any value within $[-\infty,\infty]$, it is rounded and then denoted as $\texttt{int}(z)$. 
The probability of this rounded $z$ to take to an integer value of \( k \) 
is given by:
\begin{align*}
\displaystyle \prob \left\{ \texttt{int}(z) = k \right\} & = \prob\left\{k - 0.5 \leq z < k + 0.5\right\} \\
\displaystyle & = \Phi\left(\frac{k+0.5}{\sigma}\right) - \Phi\left(\frac{k-0.5}{\sigma}\right),
\end{align*}
where \( \Phi(x) \) is the cumulative distribution function (CDF) of the standard normal distribution,
\( \displaystyle \Phi(x) = \frac{1}{2} \left[ 1 + \mathrm{erf}\left( \frac{x}{\sqrt{2}} \right) \right],
\)
with \( \mathrm{erf}(x) \) being the error function defined as
    \( \displaystyle 
    \mathrm{erf}(x) = \frac{2}{\sqrt{\pi}} \int_0^x e^{-t^2} \, \textrm{d}t
    \). 
 The explicit probability thus becomes
\begin{equation} \label{eq:exactTruncatedNormal}
\boxed{
 \displaystyle \prob \left\{ \texttt{int}(z) = k \right\} = 
 \frac{1}{2} \left[ 
\mathrm{erf}\left( \frac{k+0.5}{\sqrt{2}\sigma} \right) 
- 
\mathrm{erf}\left( \frac{k-0.5}{\sqrt{2}\sigma} \right)
\right].~ }
\end{equation}
This term can be approximated using the Mean Value Theorem for the error function, or simply by taking the first-order Taylor expansion (denote the function's argument as $x_{\pm}:= \frac{k\pm0.5}{\sqrt{2}\sigma}$): 
\[
\mathrm{erf}\left(x_{+}\right) - \mathrm{erf}\left(x_{-}\right) \approx \mathrm{erf}^{\prime}\left(\bar{x}\right) \cdot \left(x_{+}-x_{-} \right) 
\approx \frac{2}{\sqrt{\pi}} \exp{\left( -\bar{x}^2 \right)}
\]
where $\bar{x}=\frac{k}{\sqrt{2}\sigma}$ serves as the mean value.\\
This approximation holds for small $\bar{x}$ values, which translates to moderate to large standard deviations of the normal distribution (the regime where $\sigma \gtrsim 0.75$). 
Otherwise, the regime of smaller $\sigma$ values induces a probability of nearly one for no mutations, i.e., 
\( \sigma < 0.75 \leadsto \prob \left\{ \texttt{int}(z) = 0 \right\} \approx 1.0 \). 
\footnote{Another decent approximation for $\mathrm{erf}(x)$ is the \textit{hyperbolic tangent} function,
\[
\mathrm{erf}(x) \approx \mathrm{tgh}(x):= \frac{\exp(x) - \exp(-x)}{\exp(x) + \exp(-x)} ,
\]
by which the difference can be presented as
\begin{align*}
\mathrm{erf}(x_{+}) - \mathrm{erf}(x_{-}) \approx & \frac{2\left( \exp\{x_{+}-x_{-}\} - \exp\{-x_{+}+x_{-}\}\right) }{\left( \exp \{x_{+}\} + \exp\{-x_{+}\}\right) \left( \exp\{ x_{-}\} + \exp\{-x_{-}\}\right)}\\
= & \frac{2\left( \exp\{\frac{1}{\sqrt{2}\sigma} \} - \exp\{ \frac{-1}{\sqrt{2}\sigma} \}\right) }{\exp{\frac{\sqrt{2} k}{\sigma}} + \exp{\frac{1}{\sqrt{2}\sigma}} + \exp{\frac{-1}{\sqrt{2}\sigma}} + \exp{\frac{-\sqrt{2} k}{\sigma}} }
\end{align*}
without leading to a concise representation yet yielding similar numerical patterns.
}
We are now in a position to \textbf{approximate} the probability function of \eqref{eq:exactTruncatedNormal} via the following:
\begin{equation} \label{eq:approxTruncatedNormal}
\boxed{
 \displaystyle \prob \left\{ \texttt{int}(z) = k \right\} \approx
 \frac{1}{\sqrt{\pi}} \exp \left( -\frac{k^2}{2\sigma^2} \right) .}
\end{equation}
Next, we would like to express the expected $\ell_1$-norm of this TN distribution. 
Attempts to calculate the expected $\ell_1$-norm of the TN distribution in a compact form have been \textit{unsuccessful}, even upon using the approximated density function, 
\begin{align}
    \displaystyle S_{TN} := n \cdot \E \left[ |\texttt{int}(z)| \right] \approx \displaystyle \frac{n}{\sqrt{\pi}} \sum_{k=-\infty}^{+\infty} |k| \cdot \exp \left(-\frac{k^2}{2\sigma^2} \right).
\end{align}
We hypothesize that the expectation value of the continuous Normal distribution serves as a fine approximation for the discrete counterpart. 
Given the expected $\ell_1$-norm as calculated exactly for the standard Normal distribution with zero-mean,
\begin{align}
    \displaystyle \E \left[ \| \vec{z} \|_1\right] = n\cdot \int_{-\infty}^{\infty} |z| \cdot \text{pdf}(z) ~\text{d}z = n \sigma \sqrt{\frac{2}{\pi}} 
\end{align}
we empirically validated it for the TN distribution over a large range of $\sigma$ values across various dimensions. Figure \ref{fig:S-vs-sigma_TN} presents these empirical validation tests for $n=\{1,10,30,80\}$, where evidently this approximation's accuracy improves as the dimensionality increases.
We conclude that this relation, i.e., 
\begin{align}\label{eq:TN_S-sigma}
    \displaystyle S_{TN} \approx n \sqrt{\frac{2}{\pi}} \cdot \sigma_{TN} \quad \Longleftrightarrow \quad \sigma_{TN} \approx \sqrt{\frac{\pi}{2}} \cdot S_{TN}/n,
\end{align}
serves as a fine approximation per the TN distribution for the expected $\ell_1$-norm, even for small $\sigma$ values.

Finally, we mention the entropy function of this distribution using the approximated probability function of Eq.~\eqref{eq:approxTruncatedNormal}:
\begin{align}\label{eq:H_TN}
    \displaystyle H_{TN}= -\sum_{k=-\infty}^\infty  \frac{1}{\sqrt{\pi}} \exp \left( - \frac{k^2}{2\sigma^2} \right) \cdot \log_2 \left( \frac{1}{\sqrt{\pi}} \exp \left( - \frac{k^2}{2\sigma^2} \right) \right).
\end{align}

\subsubsection{The Double Geometric (DG) Distribution}
Rudolph's study proposed to use the doubly-geometric mutation, which possesses a symmetric distribution with respect to 0. 
It is achieved by drawing two random variables, $\left\{ g_1,g_2\right\}$, according to the geometric distribution, $\prob\left\{ g_{\jmath} = k \right\} = p\cdot \left(1-p \right)^k$ $(\jmath=1,2)$, and taking their difference, $z=g_1-g_2$.
The probability function of $z$ reads,
\begin{equation} \label{eq:probDGeometric}
\boxed{
\displaystyle \prob\left\{ z = k \right\} = \frac{p}{2-p} \cdot \left(1-p \right)^{\left|k\right|},\quad k \in \mathbb{Z} }
\end{equation}
with $\E \left[z\right]=0$ and $\textrm{VAR}[z]=2(1-p)/p^2$. 
Next, to generalize to the multivariate case of an $n$-dimensional mutation vector $\vec{z}$, Rudolph showed that by taking $n$ stochastically independent random variables, following the doubly-geometric distribution, the properties of symmetry and maximal entropy are kept. 
\textbf{In practice}, each random variable is drawn by the following calculation ($g_{\jmath}$ are \textit{geometrically} distributed random variables, both with parameter $p$):
\begin{equation}\label{eq:doubleGeometric}
    \begin{array}{l}
         \displaystyle g_{\jmath} \longleftarrow \left\lfloor \frac{\textrm{log} \left(1- \mathcal{U}\left(0,1\right)\right)}{\textrm{log} \left(1- p\right)}\right\rfloor \quad \jmath=1,2\\
         \displaystyle \mathcal{G}_{n}\left(0,p\right):=  g_1 - g_2. 
    \end{array}
\end{equation}
Importantly, this multivariate distribution could be \textit{controlled by the mean step-size} $S=\E \left[ \|\vec{z} \|_1\right]=n\cdot \E \left[|z_1| \right]$,
\begin{equation}\label{eq:DGstep}
\displaystyle S_{DG}(p) = n\cdot \frac{2(1-p)}{p(2-p)} \quad \Longleftrightarrow \quad p = 1 - \frac{S_{DG}/n}{\sqrt{\left(1+\left(S_{DG}/n \right)^2 \right)} + 1}.
\end{equation}
Finally, we specify the entropy function of the DG distribution:
\begin{align}\label{eq:H_DG}
    \displaystyle H_{DG}= -\sum_{k=-\infty}^\infty \frac{p}{2-p} \cdot \left(1-p \right)^{\left|k\right|} \cdot \log_2 \left( \frac{p}{2-p} \cdot \left(1-p \right)^{\left|k\right|} \right).
\end{align}
\subsubsection{Summary}
\begin{figure*}
    \centering
    \includegraphics[width=0.45\columnwidth]{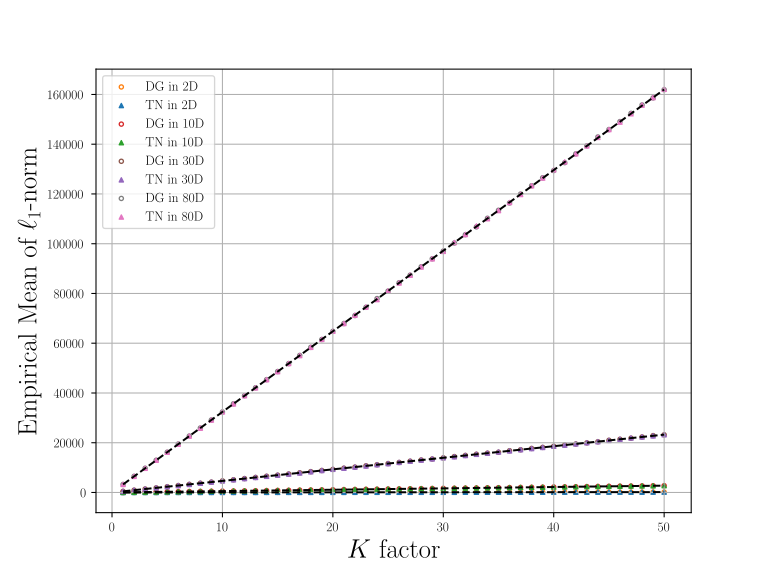} ~
    \includegraphics[width=0.45\columnwidth]{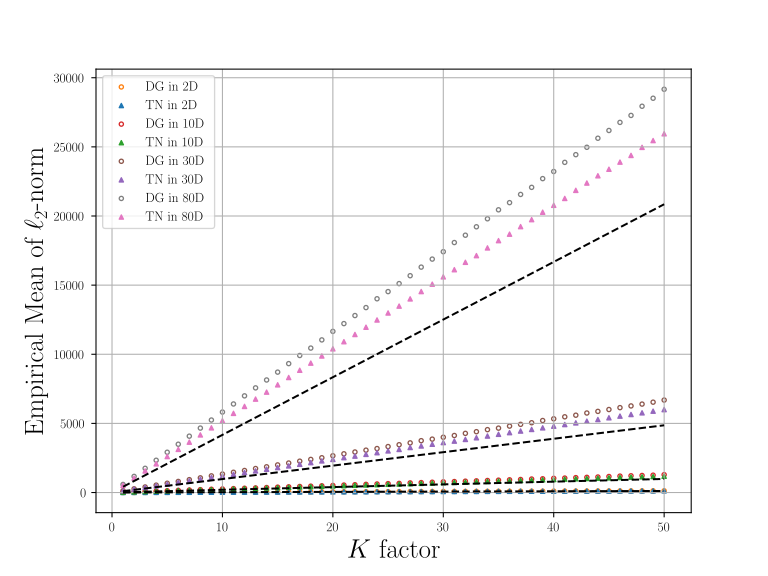}
    \caption{Assessing the DG and TN distributions: The empirical mean of the $\ell_1$- versus the $\ell_2$-norm of populations of randomly generated $n$-dimensional integer vectors adhering to an ellipsoid form of the step-sizes, that is $S_i=K\cdot i$ subject to a factor $K \in \{1,\ldots , 50 \}$ across multiple dimensions $n\in \{ 2,10,30,80\}$. The theoretical step-sizes are depicted as the black dashed lines. [LEFT] $\ell_1$-norm, and [RIGHT] $\ell_2$-norm. Evidently, the $\ell_2$-norm is completely off when $n\gtrsim 30$, being unable to quantify the vectors.}
    \label{fig:ExpectedL1norm}
\end{figure*}

We focus on the TN and DG distributions and numerically assess the accuracy of the aforementioned formulations. 
We empirically compute the mean $\ell_1$- and $\ell_2$-norm of randomly generated $n$-dimensional integer vectors $S$ per the four distributions. 
To this end, we consider an \textit{ellipsoid} form of the step-sizes and inflate it using a spectrum $K\in \{1,\ldots , 50 \}$ across multiple dimensions $n\in \{ 2,10,30,80\}$, that is
\[ 
S_i=K\cdot i \quad \forall i\in 1\ldots n~.
\]
This numerical evaluation is presented as Figure \ref{fig:ExpectedL1norm}, with the theoretical expected $\ell_1$ step-size being $S=K\cdot \sum_i i = \frac{1}{2}Kn(n+1)$.
The results corroborate the validity of our formulations per the $\ell_1$-norm, while showing a good fit between the theoretical and the empirical mean values, scalable across $K$ and $n$. At the same time, this basic numerical evaluation demonstrates the \textbf{inability of the $\ell_2$-norm to quantify the integer vectors over the discrete lattice}, since its accuracy deteriorates as the dimensionality increases ($n\gtrsim 30$).

\begin{figure*}
    \centering
    \begin{tabular}{ c  c  c}
    \hline
    $S=1$ & \includegraphics[width=0.45\columnwidth]{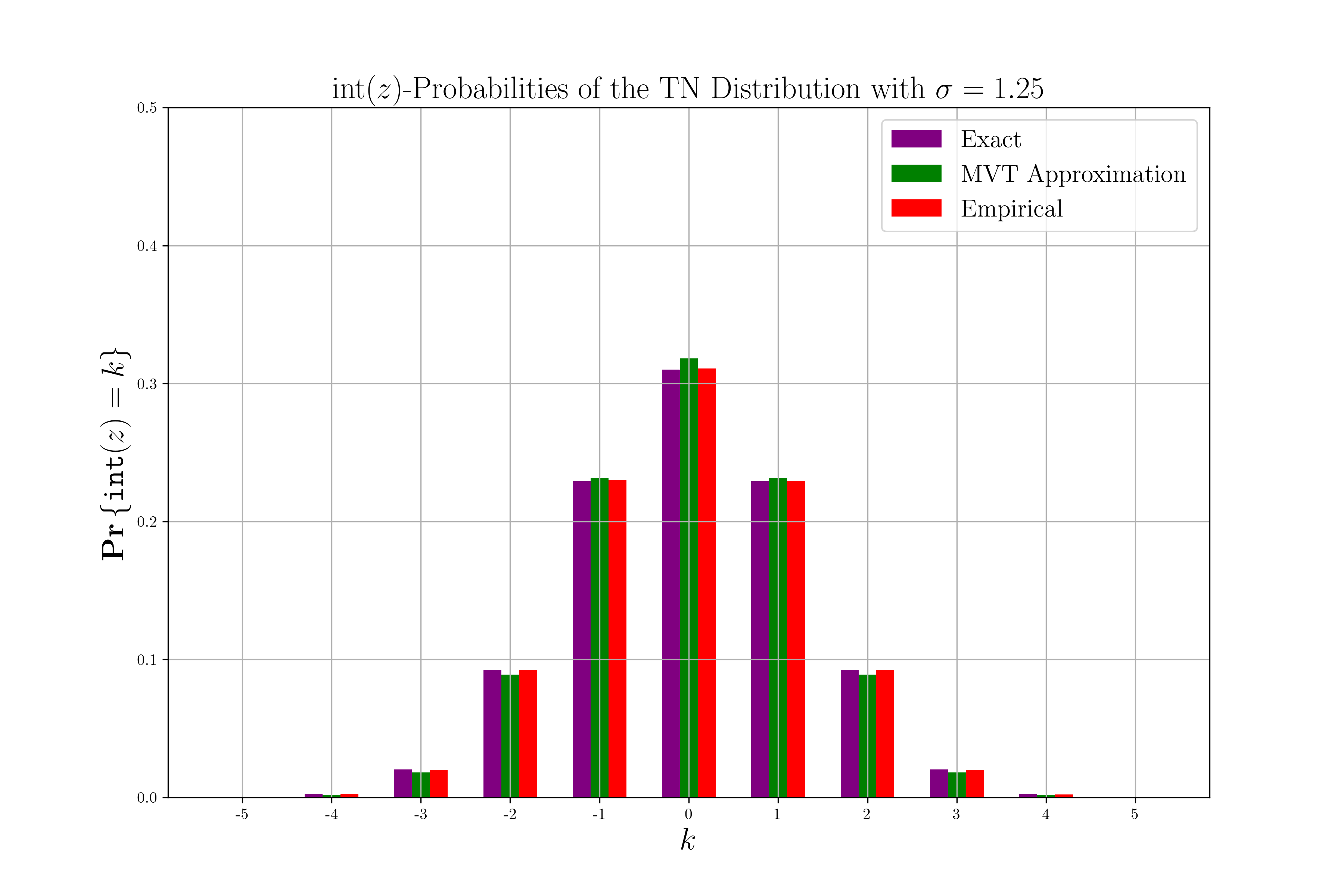} & \includegraphics[width=0.45\columnwidth]{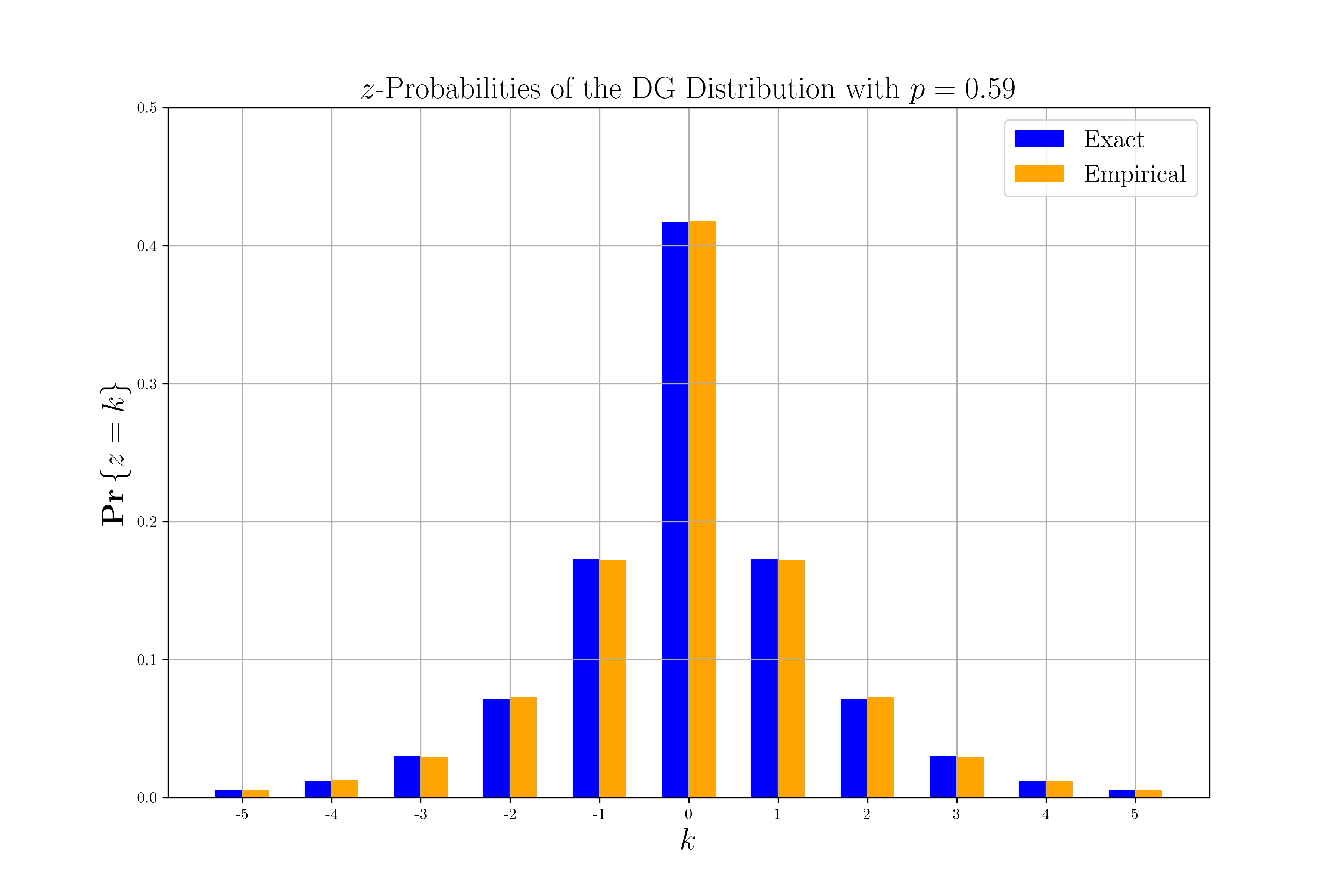}\\
    \hline 
    $S=2$ & \includegraphics[width=0.45\columnwidth]{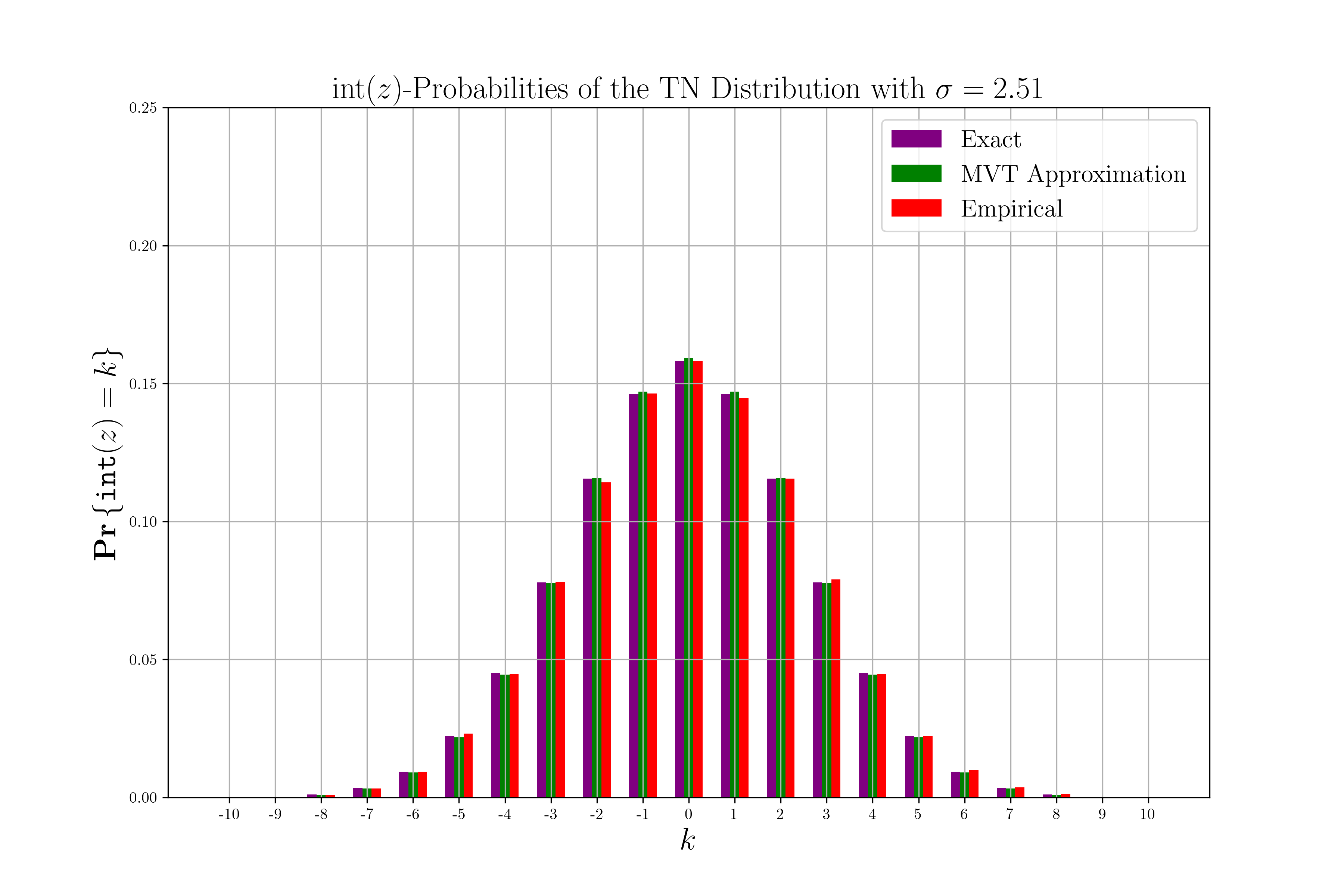} &  \includegraphics[width=0.45\columnwidth]{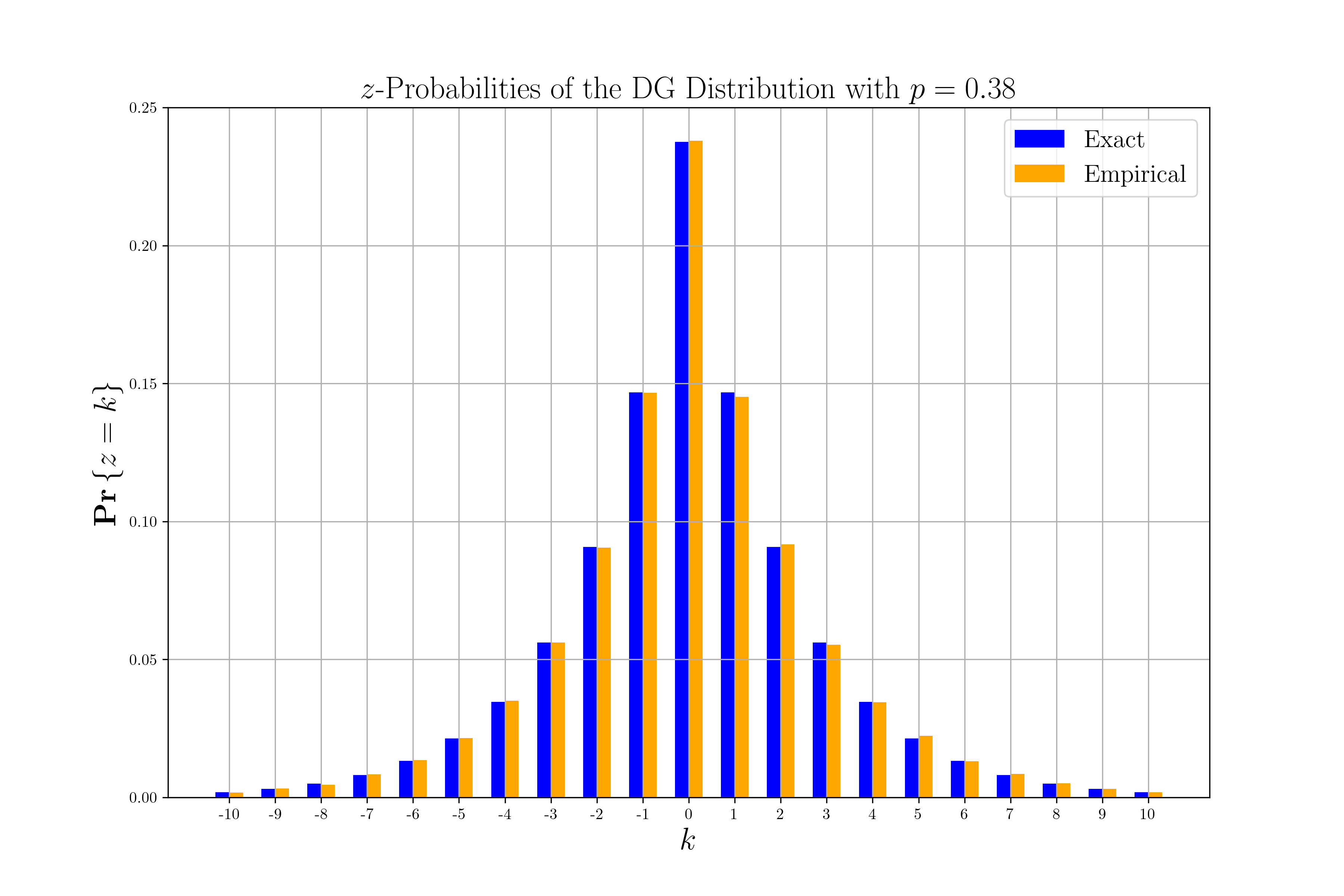}\\
    \hline 
    $S=3$ & \includegraphics[width=0.45\columnwidth]{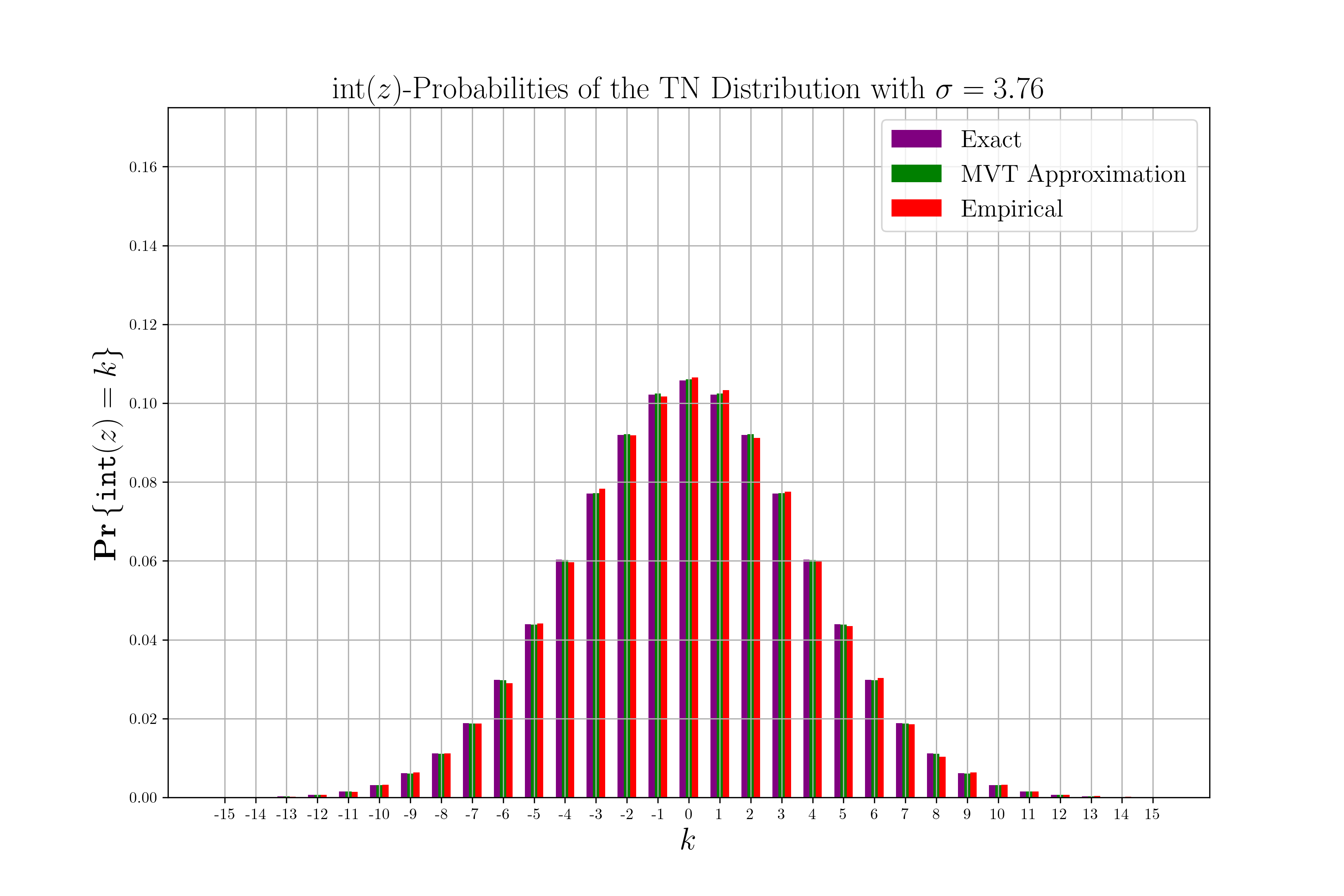} & \includegraphics[width=0.45\columnwidth]{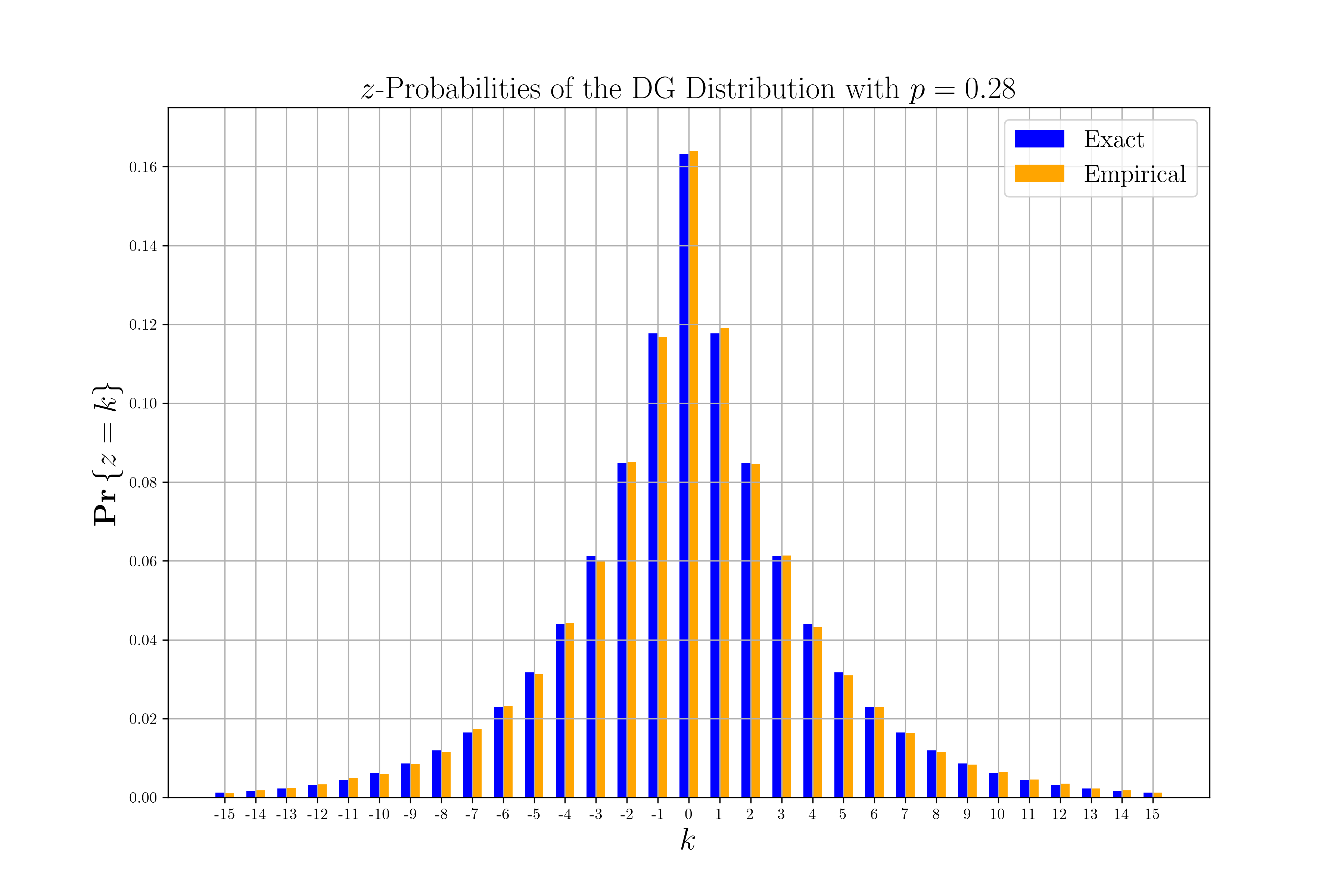} \\
    \hline
    $S=4$ & \includegraphics[width=0.45\columnwidth]{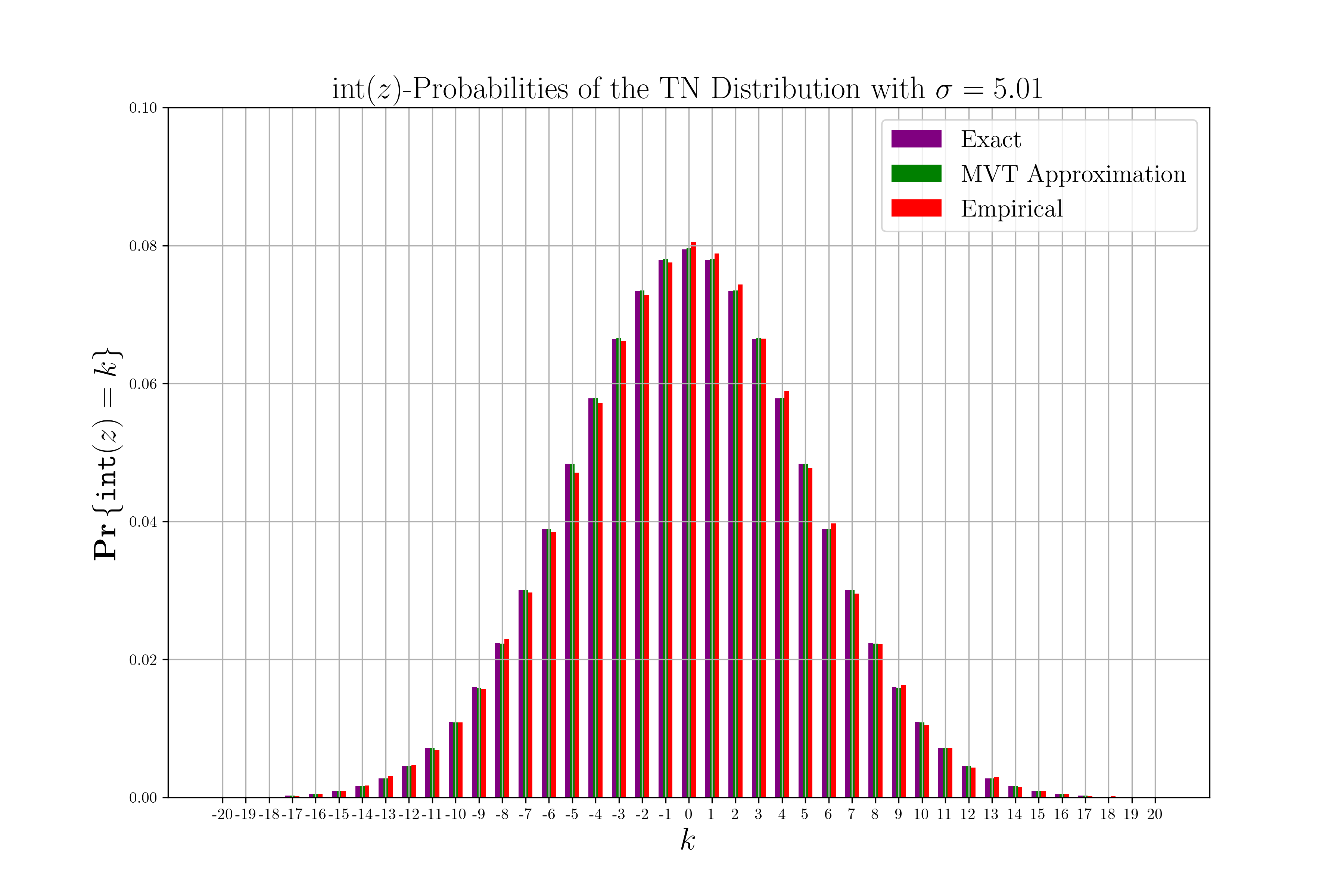} &  \includegraphics[width=0.45\columnwidth]{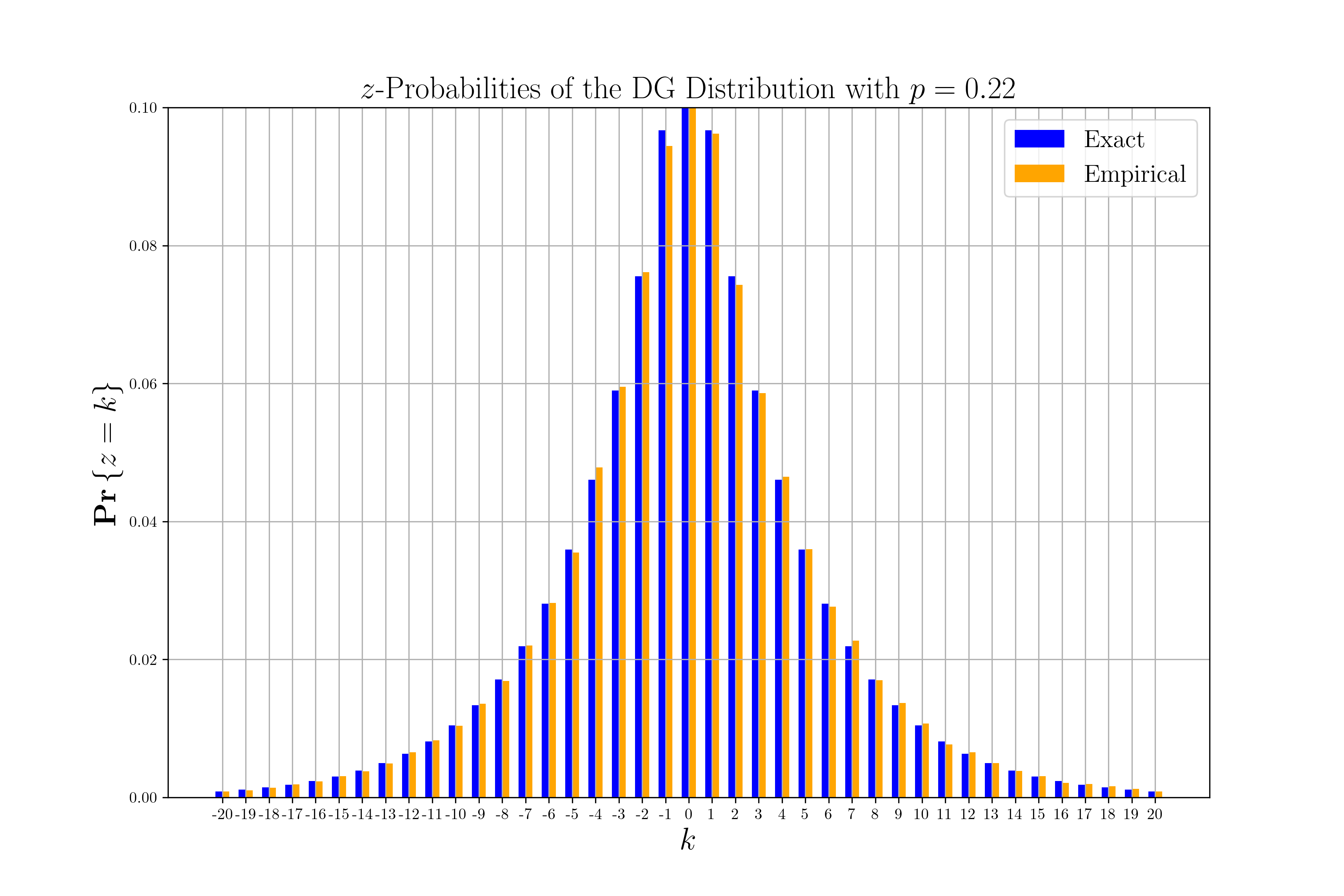} \\
    \hline 
    \end{tabular}
    \caption{The single-variable TN and DG distributions per four concrete mean step-size values $S\in\{1,2,3,4\}$ (top to bottom): [LEFT] TN probabilities (exact versus MVT approximation versus empirical) using $\sigma_{TN}(S)$ from Eq.~\eqref{eq:TN_S-sigma} (explicit $\sigma$ values are specified on each plot); [RIGHT] DG probabilities (exact versus empirical) per the computed $p$ when accounting for $S_{DG}$ using \eqref{eq:DGstep} (explicit $p$ values are specified on each plot).} \label{fig:probabilitiesValidation}
\end{figure*}

Next, by setting $S$ as the driving step-size, and by exercising the relations $\sigma_{TN}(S)$ and $p_{DG}(S)$ (using Eqs.~\eqref{eq:TN_S-sigma} and \eqref{eq:DGstep}, respectively), we compare the two distributions over a set of $S$ values and present their histograms as a gallery in Figure \ref{fig:probabilitiesValidation}. Evidently, their shapes differ -- while the TN exhibits relative smoothness and strong resemblance to the convex Bell curve as expected, the DG distribution possesses a spiked, concave shape with a pronounced peak at zero. 
%
\section{Correlated Samples via Rotations} \label{sec:correlations}
The challenge of well-defining rotations over the integer lattice stems from the fundamental incompatibility between the geometric concept of rotation, which involves trigonometric functions that induce irrational coordinate values, and the discrete, axis-aligned structure of the lattice. 
An ambition to define $\ell_1$-norm-preserving rotations is unrealistic in the generalized case, and can be fulfilled only for a limited set of special cases:
\begin{enumerate}[nosep, label={(\roman*)}]
   \item Permutation of coordinates, e.g., in the 2D case via the transformation $$R_{\text{perm}}:= 
   \begin{bmatrix}
        \hfill 0 & \hfill 1 \\ 
        \hfill 1 & \hfill 0
    \end{bmatrix}.$$   
    \item Sign changes of coordinates, e.g., in the 2D case via the transformation $$R_{\text{sign}}:= \begin{bmatrix}
        \hfill -1 & \hfill 0 \\ 
        \hfill 0 & \hfill 1
    \end{bmatrix}.$$    
    \item Combinations of the above (i)+(ii).
\end{enumerate}
Next, we will propose a rotation transformation for integer search. 
\subsection{Proposed Rotation Transformation}
Dropping the requirement for $\ell_1$-norm preservation, given an uncorrelated mutation instance $\vec{z}_u$, a correlated mutation instance $\vec{z}_c$ is attainable by \textit{rounding to the nearest integer} the sequence of $n(n-1)/2$ rotations using the operator $\mathbf{R}(\theta):=(r_{k\ell})$
\begin{align}
    \displaystyle \vec{z}_c = \texttt{round}\left[ \left(\prod_{i=1}^{n-1} \prod_{j=i+1}^{n} \mathbf{R}(\alpha_{ij})  \right) \cdot \vec{z}_u \right] ~ .
\end{align}
$\mathbf{R}$'s matrix form is identical to the unity, except for 4 elements:
\[
\displaystyle r_{kk}=r_{\ell\ell}=\cos(\alpha_{k\ell}), \quad \quad r_{k \ell}=-r_{\ell k}=-\sin(\alpha_{k\ell}),
\]
where the angle $r_{k \ell}$ represents a correlation measure between the $k^{th}$ and $\ell^{th}$ random variables. 
Overall, this procedure is well-defined when given an uncorrelated mutation vector $\vec{z}$ and a vector of rotational angles $\vec{\alpha}$:
\begin{align*}
\boxed{
\begin{array}{l}
     \textbf{\texttt{rotateInt}} \left( \vec{z},~\vec{\alpha} \right) \\
    \quad \textbf{for~} j=1,\ldots,n\cdot (n-1)/2~\textbf{~do} \\
    \quad \quad \vec{z} \longleftarrow \mathbf{R}(\alpha_j) \vec{z} \\
    \quad \textbf{end} \\
    \quad \textbf{return } \left\{ \texttt{round}\left(\vec{z} \right) \right\}
    \end{array}
    }
\end{align*}
\paragraph{The Original Continuous Usage}
A ``pure'', rounding-free transformation, i.e., \[ \vec{z}_c = \left(\prod_{i=1}^{n-1} \prod_{j=i+1}^{n} \mathbf{R}(\alpha_{ij})  \right) \cdot \vec{z}_u ,\] played the rotation role for generating correlated normal mutations in Schwefel's formulation of the Standard ES \cite{Schwefel,Baeck-book}. 
This transformation was designed to store covariance information by means of the $n$-dimensional variances' vector $\vec{\sigma}$ as well as the $n(n-1)/2$-dimensional vector of rotational angles $\vec{\alpha}$. 
In the continuous domain, the statistical covariance of the decision variables is meaningful, also due to its hypothesized relation to the inverse Hessian matrix of the search landscape (which was later confirmed \cite{Shir-Yehudayoff_TCS2020}). 
Therefore, it was important to design transformations that account for such information and store it. Given decision variables $i$ and $j$, the transformation from the covariance element $c_{ij}$ to the rotational angle $\alpha_{ij}$, with $c_{ii}\equiv \sigma_i^2$, constitutes another relation that proved useful:
\begin{align} \label{eq:tanAlpha}
\displaystyle \alpha_{ij} = \frac{1}{2} \arctan \left( \frac{2c_{ij}}{\sigma_i^2-\sigma_j^2} \right) ~ ,
\end{align}
where $\alpha_{ij}=0$ whenever no correlation exists.
Indeed, this representation of variances and angles was preferred in the Standard ES over the canonical matrix form $(c_{ij})$ since it was easier to maintain the mathematical property of the covariance being positive-definite when the update occurs via self-adaptation.
Later on, Rudolph \cite{Rudolph92} verified the validity of this representation and showed that the usage of such angles does not restrict the generality of correlated ES mutations.

In our context of integer search, we can use the aforementioned \texttt{rotateInt} procedure to generate correlated integer mutations. 
However, our usage is not concerned by any mathematical constraint to form a proper covariance matrix, whose interpretation is limited and usefulness unknown over the integer lattice. \\
Such a procedure is proposed as Algorithm \ref{algo:mutate}, relying on the \texttt{genUncorrelatedMutation()} subroutine, which generates the uncorrelated samples according to the type (denoted as \{\texttt{DG}, \texttt{TN}\}, although a single call of \texttt{DG} obtains a Geometric sample and hence the need to call it twice and take the difference).\\ 
We will numerically evaluate this procedure in the next subsection.\\
\IncMargin{0.5em}
{\LinesNotNumbered
\begin{algorithm}[h]
\caption{A procedure to generate a correlated mutation vector of integers. $\left\{\vec{S},\vec{\alpha}\right\}$ are the step-sizes and rotation angles, respectively. 
This procedure handles either the DG or the TN distributions via its subroutine \texttt{genUncorrelatedMutation()}, which accounts for the type (denoted as \{\texttt{DG}, \texttt{TN}\}). The rotation, via \texttt{rotateInt()}, is agnostic with respect to the underlying distribution of its uncorrelated input sample. 
\label{algo:mutate}}
\normalsize 
\Indm\SetKwFunction{Func}{\textbf{corrMutate}}
\Func{$\vec{S},~\vec{\alpha},~n,~\texttt{type}$}\;
\Indp
    $\vec{z}_u \longleftarrow  \texttt{genUncorrelatedMutation}\left( \vec{S},~\texttt{type} \right)$ \;
    $\vec{z} \longleftarrow \texttt{rotateInt}\left( \vec{z}_u,~\vec{\alpha}\right) $ \;
    \If{\texttt{type==DG}} { 
        $\vec{z}_g \longleftarrow \texttt{genUncorrelatedMutation}\left( \vec{S},~\texttt{type} \right)$ \;
        $\vec{z}_g^{\prime} \longleftarrow \texttt{rotateInt}\left( \vec{z}_g,~\vec{\alpha}\right) $ \;
        $\vec{z} \longleftarrow \vec{z} - \vec{z}_g^{\prime}$ // difference of two geometric samples
    }  
\Return {$\left\{\vec{z}\right\}$}
\end{algorithm}}
\DecMargin{0.5em}
Notably, it is possible to artificially form a matrix using a step-size vector $\vec{S}$ and rotational angles $\vec{\alpha}$ by exercising Eqs.~\eqref{eq:TN_S-sigma} and \eqref{eq:tanAlpha}, respectively. However, this construction is not guaranteed to hold the necessary covariance matrix mathematical properties, with positive semidefiniteness (PSD) being the primary.\\ 
\textbf{Whenever it is mathematically sound}, after applying a numerical procedure to enforce the PSD property, we will consider such a constructed covariance matrix and denote it by $\mathbf{C}_{FTN}$ (which stands for Forced TN). 
Then, we can use it to generate correlated multivariate normally distributed vectors with truncation:
\begin{align}
    \displaystyle \vec{z}_{FTN} &  \sim\texttt{round} \left[ \mathcal{N}\left(\vec{0}, \mathbf{C}_{FTN}\right) \right]~.
\end{align}
\subsection{Numerical Assessment of 2D Rotations}
\begin{figure*}[h]
    \centering
    \includegraphics[width=1\columnwidth]{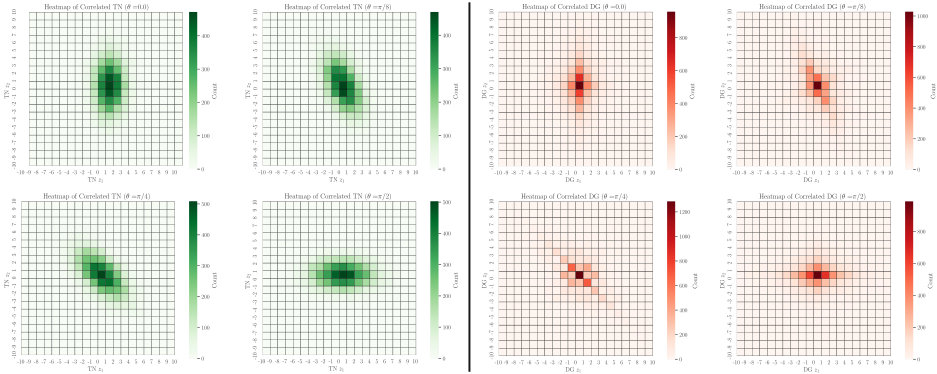}
    \caption{Heatmaps depicting populations of rotated 2D samples with $S_1=1,~S_2=2$ and $\theta \in \{0,\frac{\pi}{8},\frac{\pi}{4},\frac{\pi}{2} \}$ (set in this order clockwise - see titles). [LEFT, green]: TN distribution, [RIGHT, red]: DG distribution. 
    The population size was set to $10^4$.}
    \label{fig:rotated2D}
\end{figure*}
\paragraph{Heatmaps} To visually demonstrate the ability to generate correlated integer mutations, we present populations of 2D samples per the TN and DG distributions across multiple rotation angles with $S_1=1,~S_2=2$ and $\theta \in \{0,\frac{\pi}{8},\frac{\pi}{4},\frac{\pi}{2} \}$ -- via heatmaps (see Figure \ref{fig:rotated2D}). 
The DG distribution exhibits a spiked pattern, in contrast to the smooth TN pattern, where mind should be given to the scale. 
This pattern constitutes a two-dimensional extension of the single variable case, which also demonstrated discrete versus smooth profiles, respectively (Figure \ref{fig:probabilitiesValidation}).
\begin{figure}
    \centering
    \includegraphics[width=0.7\columnwidth]{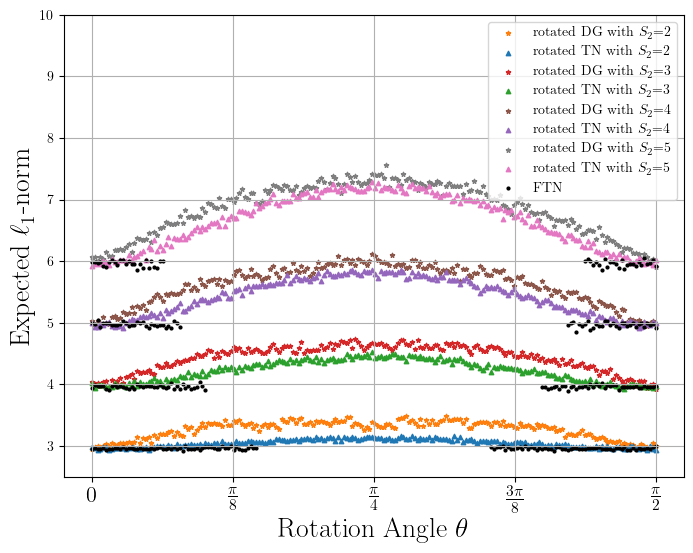}
    \caption{The empirical $\ell_1$-norm of randomly generated 2D integer vectors per the TN/DG distributions when subject to systematic rotation with $\theta \in [0,\frac{\pi}{2}]$ under four settings $S_2 \in \{2,3,4,5\}$ with $S_1=1$. 
    The population size was set to $10^4$.}
    \label{fig:thetaRotatedL1}
\end{figure}

\paragraph{$\ell_1$-norm under systematic rotations} 
Next, we randomly generated populations of 2D integer vectors using the two distributions when subject to \textit{systematic} rotation by continuous angle $\theta \in [0,\frac{\pi}{2}]$ under four settings $S_2 \in \{2,3,4,5\}$ with $S_1=1$. We empirically compute the mean $\ell_1$-norm of such rotations in Figure \ref{fig:thetaRotatedL1}, where the FTN samples are partially defined over this range of $\theta$.
It is evident that the proposed rotation does not preserve the $\ell_1$-norm, which is to be expected, and that the bias is maximized at the $\theta=\frac{\pi}{4}$ rotational angle. 
Also, the DG distribution consistently obtains higher $\ell_1$-norm values with a less smooth pattern when compared to the TN distribution. 

\begin{figure*}
\centering
    \includegraphics[width=0.45\columnwidth]{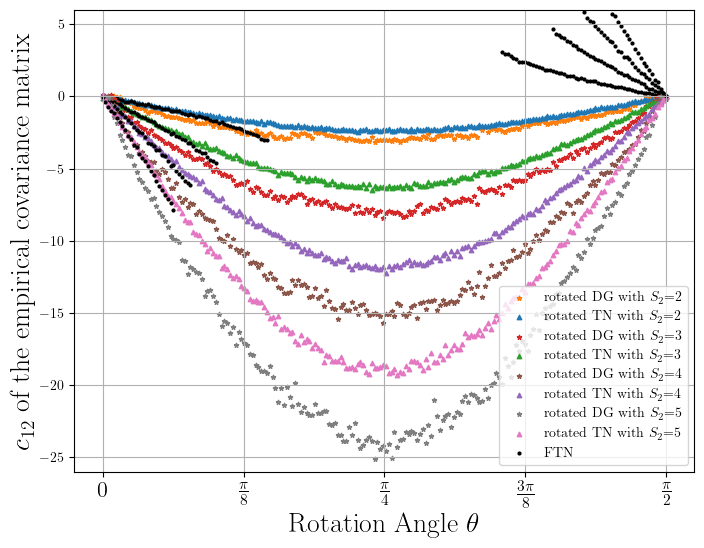} ~~
    \includegraphics[width=0.45\columnwidth]{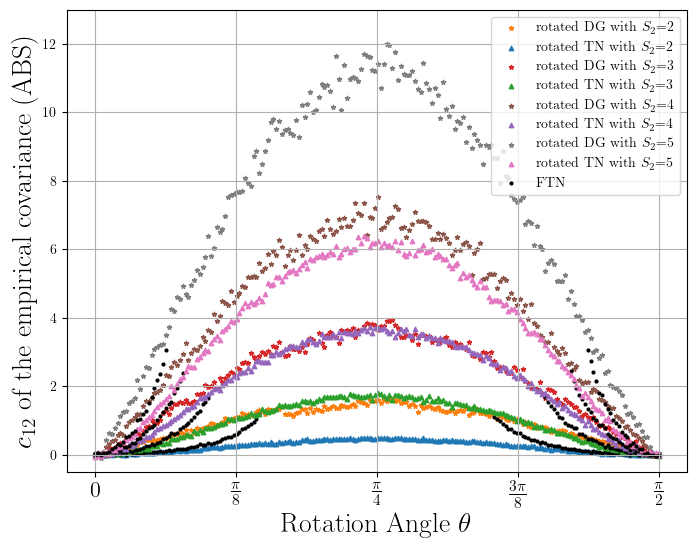}
    \caption{Numerical assessment of various correlation measures over randomly generated 2D integer vectors per the TN/DG distributions when subject to systematic rotation with $\theta \in [0,\frac{\pi}{2}]$ under four settings $S_2 \in \{2,3,4,5\}$ with $S_1=1$:\\ (LEFT) default statistical covariance, and (RIGHT) covariance over absolute values ($|z_i|$). The FTN samples (depicted by black points) are defined only partially over this range. }
    \label{fig:thetaRotatedcorr}
\end{figure*}

\paragraph{Statistical measures for correlations} 
Finally, we tested four statistical measures to quantify the correlation of the two integer decision variables in the 2D case. 
We generated the populations similarly to the previous simulation of the $\ell_1$-norm under rotations. 
The Pearson correlation measure failed to provide a fair representation of the data and is excluded from this report. 
Figure \ref{fig:thetaRotatedcorr} presents two correlation measures over randomly generated 2D integer vectors per the TN/DG distributions when subject to systematic rotation by $\theta \in [0,\frac{\pi}{2}]$ under four settings $S_2 \in \{2,3,4,5\}$ with $S_1=1$: the default statistical covariance, and the covariance over absolute values (that is, empirical covariance over $|z_i|$ values). We also computed $\ell_1$-driven correlation via the LASSO coefficients \cite{LASSO-Tibshirani1996} (see Supplementary Material). 
The FTN samples (depicted by black points in both figures) are defined only partially over this range of $\theta$. 
The DG distribution consistently exhibits abrupt changes in its calculated measures, but its trend curves are consistent. 

Setting aside the FTN data points (which are anyway not well-defined), it is evident that the three statistical measures are able to quantify the correlation proportionally to the rotation in an expected pattern. 
Assessing the effectiveness of the considered statistical measures in meaningfully expressing correlations driven by the $\ell_1$-norm requires further research and is not addressed here.

\section{Entropy} \label{sec:entropy}
Rudolph investigated the family of maximum entropy \textbf{discrete} distributions \cite{rudolph1994evolutionary} by addressing the following optimization problem:
\begin{align}\label{eq:maxEntropy}
\boxed{
\begin{array}{ll}
\medskip 
\displaystyle \textrm{maximize}_{p_k} & \displaystyle H:= -\sum_{k=-\infty}^\infty p_k \log p_k\\
\displaystyle \textrm{subject to: }& \displaystyle p_k = p_{-k}, \quad \forall k \in \mathbb{Z},\\ 
\displaystyle & \displaystyle \sum_{k=-\infty}^\infty p_k = 1,\\
\displaystyle & \displaystyle \sum_{k=-\infty}^\infty k^2 p_k = \sigma^2,\\ 
\displaystyle & \displaystyle p_k \geq 0, \quad \forall k \in \mathbb{Z}.
\end{array} }
\end{align}
By treating the problem's \textbf{Lagrangian} (neglecting the last condition)
\begin{equation}
\displaystyle \mathcal{L}(p, a, b) = -\sum_{k=-\infty}^\infty p_k \log p_k + a \left( \sum_{k=-\infty}^\infty p_k - 1 \right) + b \left( \sum_{k=-\infty}^\infty |k|^m p_k - \sigma^2 \right)
\end{equation}
Rudolph concluded that the DG distribution constitutes the \textit{maximizer} of this optimization problem.
Next, we would like to assess the entropy functions of the considered distributions. 
\subsection{Exact Entropy versus Mean Step-Size (1D)}
\begin{figure}
    \centering
    \includegraphics[width=0.75\columnwidth]{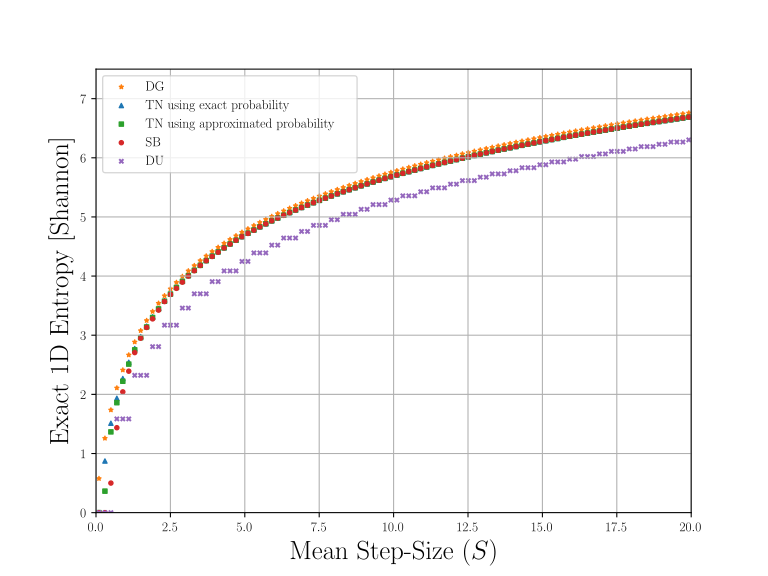}
    \caption{The entropy function of the single-variable distributions over the spectrum of $S$, which controls their defining step-size: $H_{TN}$ \eqref{eq:H_TN} via $\sigma_{TN}(S)$ from Eq.~\eqref{eq:TN_S-sigma}, using the exact \eqref{eq:exactTruncatedNormal} and the approximated \eqref{eq:approxTruncatedNormal} probabilities, which yield strong resemblance; $H_{DG}$ \eqref{eq:H_DG} over $p$ (when accounting for $p\in[0,1]$ using the relation $S_{DG}(p)$ from \eqref{eq:doubleGeometric}), and additionally $H_{SB}$, using its exact probability function \eqref{eq:shiftedBinomial} and by exercising \eqref{eq:NbyS-Bimonial} -- showing a strong resemblance to the TN distribution, as expected from the Binomial-Normal relation.}
    \label{fig:entropy}
\end{figure}

In order to establish a baseline for comparing the distributions, a common mutation step-size must be well-defined. 
We are able to draw a comparison of the entropy as a function of the integer mutation step-size, by considering the mean step-size $S$ as the driving parameter to dictate $\sigma_{TN}$ and $p$. 
Accordingly, we compare the entropy functions that are derived from the two distributions -- $H_{TN}$ and $H_{DG}$ -- 
and additionally consider $H_{SB}$ (using its exact probability function \eqref{eq:shiftedBinomial}) by exercising \eqref{eq:NbyS-Bimonial} as well as $H_{DU}$ \eqref{eq:HDU}. \\ 
Figure \ref{fig:entropy} depicts the entropy functions over a spectrum of $S$, using the specified relations. 
While the DU distribution exhibits the lowest entropy values, the entropy functions of the SB and TN distributions practically merge. Most importantly, it is evident that the DG distribution consistently possesses the highest entropy values, constituting a numerical validation to Rudolph's result.
\subsection{Estimated Entropy of 2D Samples}
\begin{figure}
    \centering
    \includegraphics[width=0.75\columnwidth]{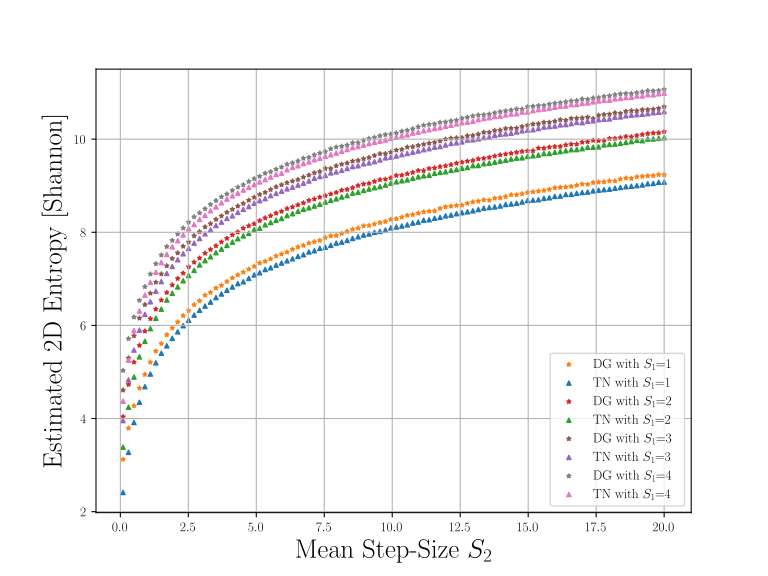} 
    \includegraphics[width=0.75\columnwidth]{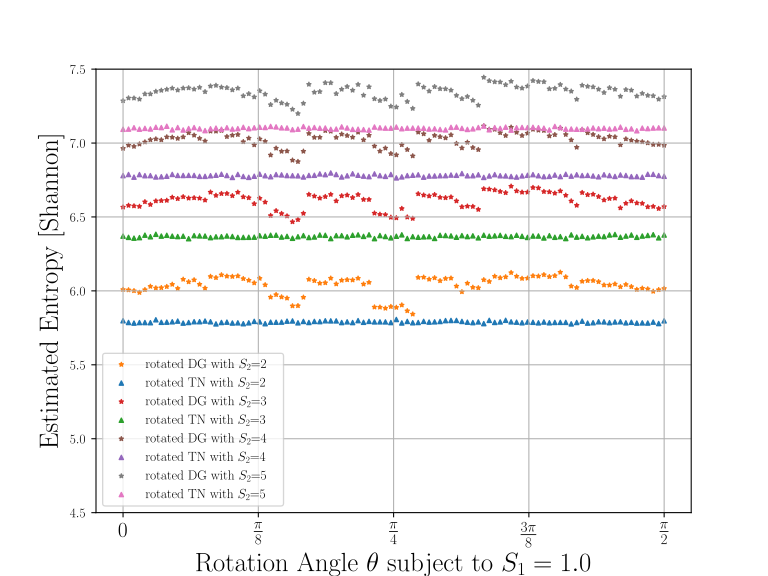} 
    \caption{Entropy histogram-based estimations of the 2D DG and TN distributions when controlled by $S_1$ and $S_2$. 
    Uncorrelated samples with multiple step-sizes: $S_1\in \{1,2,3,4 \}$ over a spectrum of $S_2$ values (top), and correlated samples subject to $S_1=1$ and $S_2\in\{2,3,4,5 \}$ over a spectrum of rotation angles $\theta \in [0,\pi/2]$ (bottom). Evidently, the DG distribution consistently possesses higher entropy values in all the 2D scenarios, both correlated and uncorrelated.}
    \label{fig:entropy2D}
\end{figure}

In order to numerically assess the entropy values of 2D randomly generated integer vectors, we turn to histogram-based estimations of the generated populations.
We consider 2D populations drawn from the DG and TN distributions when controlled by $S_1$ and $S_2$. 
Firstly, uncorrelated samples are generated with $S_1=1$ over a spectrum of $S_2$ values, and their estimated entropy values are presented in Figure \ref{fig:entropy2D}[top].  
Secondly, correlated samples subject to $S_1=1$ and $S_2\in\{2,3,4 \}$ are generated over a spectrum of rotation angles $\theta \in [0,\pi/2]$ and undergo histogram estimation to obtain the entropy values -- depicted in Figure \ref{fig:entropy2D}[bottom].  

\textbf{Evidently, the samples generated by the DG distribution consistently possess higher entropy values in all the investigated 2D scenarios -- suggesting that Rudolph's result holds also under the proposed rotation transformation. }


%
%
\section{Practical Observations}\label{sec:experiments}
We present empirical observations of IES variants across quadratic models and discuss them in light of the aforementioned theoretical principles. 
These observations are based on a recent empirical study that investigated the performance of IESs over integer quadratic programming (IQP) problems \cite{Shir-Emmerich_gecco2025}.
\subsection{Standard IES}
By adopting the proposed correlated mutation operator (\texttt{corrMutate} of Algorithm \ref{algo:mutate}), it is necessary to define the self-adaptation scheme of the strategy parameters $\vec{S}$ and $\vec{\alpha}$. 
Adhering to Schwefel's original scheme, the update steps read
\begin{align}
& S_i^{\prime}\longleftarrow S_i\cdot \exp\left\{ \tau_g\cdot\mathcal{N}_g + \tau_{\ell}\cdot\mathcal{N}_i\left(0,1\right)\right\}  \quad  \text{for } i=1\ldots n~, \\
& \alpha_j^{\prime}\longleftarrow \alpha_j + \beta \cdot \mathcal{N}_{j}\left(0,1\right) \quad  \text{for } j=1,\ldots,n\cdot (n-1)/2~,
\end{align}
using the default parametric settings: 
\( \tau_g:=\frac{1}{\sqrt{2\cdot n}},~\tau_{\ell}:=\frac{1}{\sqrt{2\cdot \sqrt{n}}}\) 
and $\beta:= 0.0873$ \cite{Baeck-book}. 
Additionally, the renowned 1/5th success-rule \cite{Baeck-book} for the step-size adaptation is considered in play with either the TN or DG mutation distributions.
Altogether, six IESs variants were considered, reflecting the combinations across the underlying distributions DG/TN and the application of rotations: \text{(1+1)-DG}, \text{(1+1)-TN}, \text{correlatedDG}, \text{correlatedTN}, \text{uncorrelatedDG}, and \text{uncorrelatedTN}. 
The DU distribution has been tested and consistently failed to deliver competitive results (numerical results are excluded), implying that the distribution plays an important role.\\
Analyses have been carried out by the \texttt{IOHanalyzer} tool \cite{IOHanalyzer}.

\subsection{Preliminary: Stagnation over the Sphere}
\begin{figure}
    \centering
    \includegraphics[width=0.9\columnwidth]{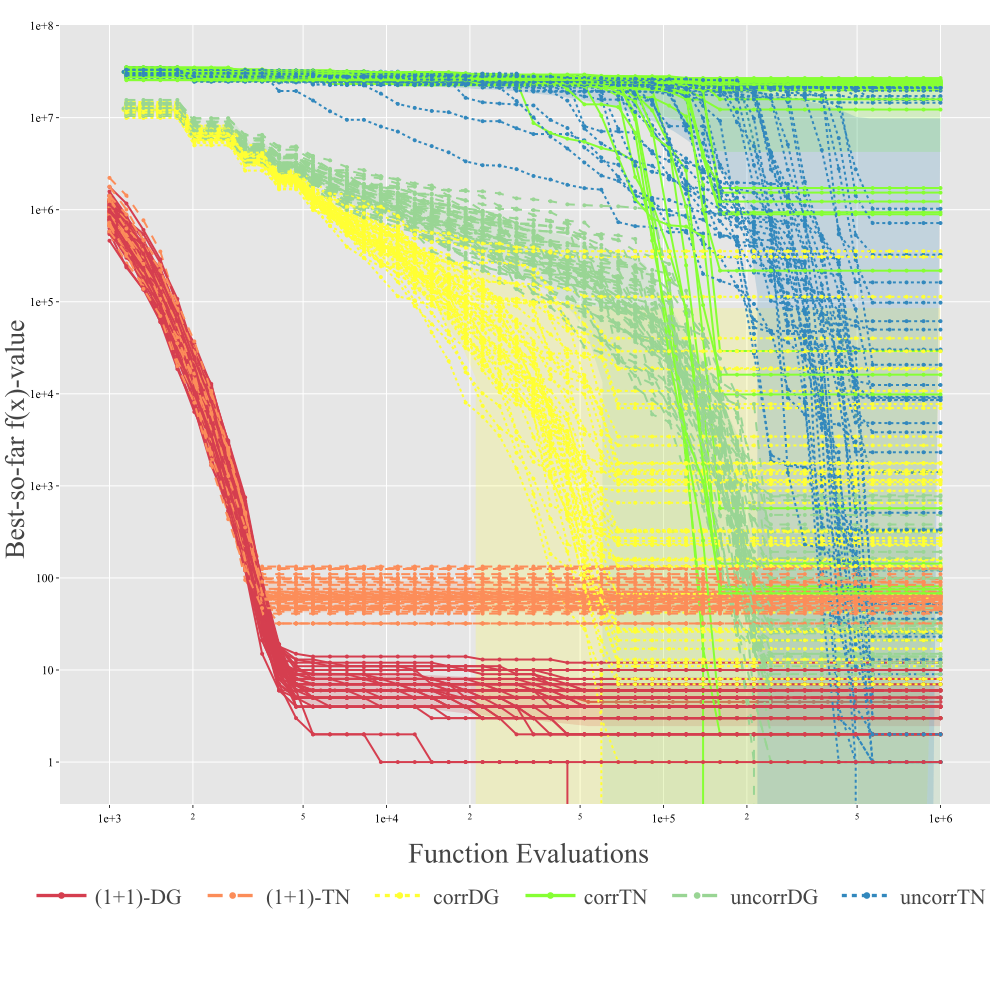}
    \caption{The six IESs variants applied to the 80D integer Sphere function: the Expected Target Value plot of the fixed-budget perspective over the individual 50 runs per strategy.}
    \label{fig:80Dfv}
\end{figure}

The six IESs were tested over the integer Sphere function by means of 50 runs over three dimensionalities $n\in \left\{30,50,80 \right\}$. 
The summary of the 80D runs via the fixed-budget perspective is provided in Figure \ref{fig:80Dfv} as an Expected Target Value plot \cite{IOHanalyzer}.
The (1+1)-DG clearly outperforms the other strategies on this test-case, but it is also evident that all six IESs do not reach the optimum in most of the runs. Instead, they enter stagnation phases away from the optimum, even when no local optima exist. 
This behavior is consistently observed across other dimensionalities. 
Thus, we argue that this stagnation stems from inherent issues within the IESs' underlying mechanisms, independent of the search landscape's difficulty. 

\subsection{Integer Quadratic Programming (IQP)}
\begin{figure}[h]
    \includegraphics[width=0.5\columnwidth]{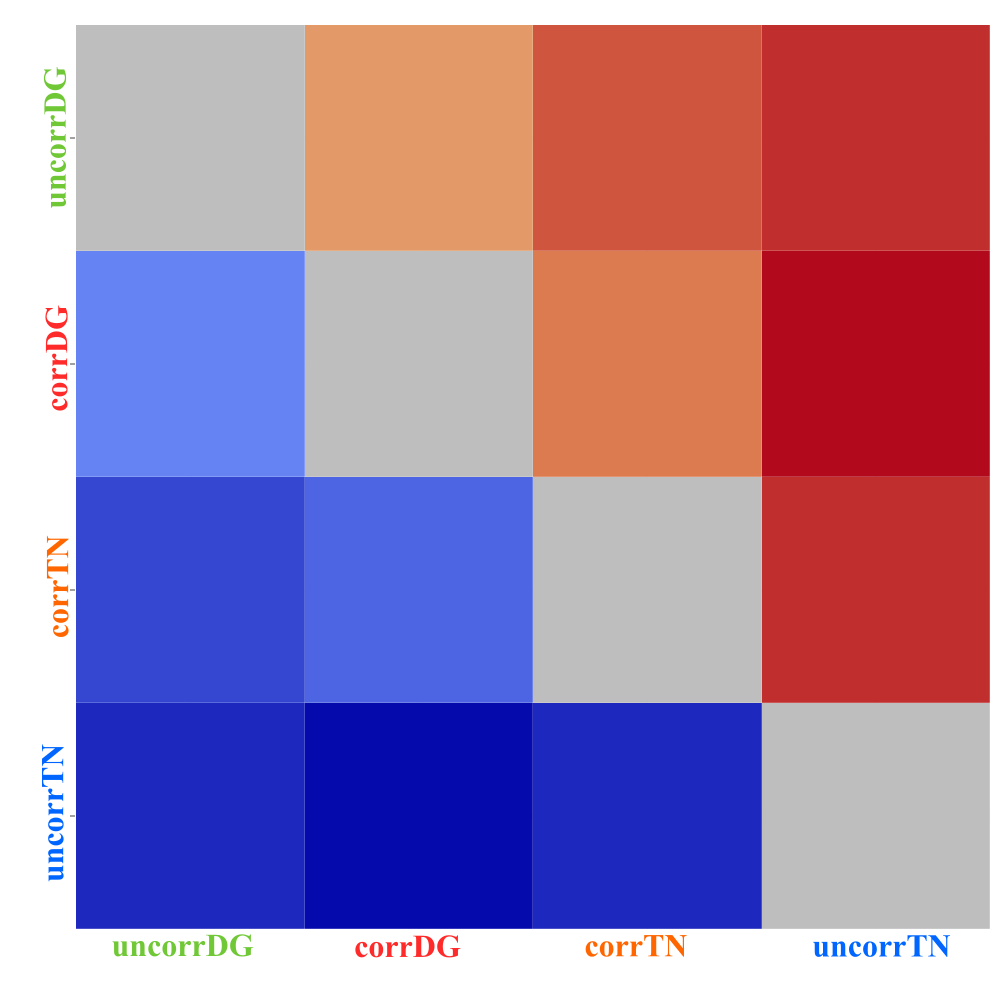} ~
    \includegraphics[width=0.5\columnwidth]{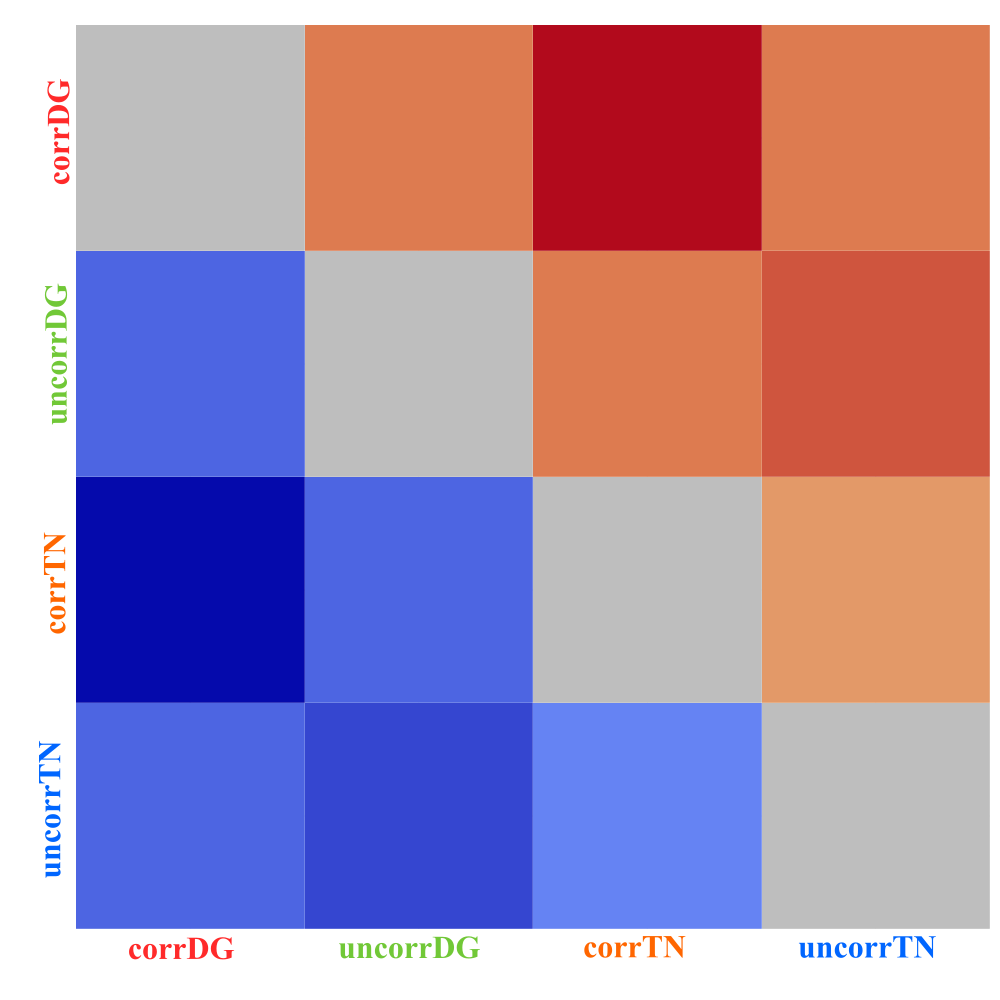}
    \caption{IESs applied to IQP: pairwise comparisons among the four IESs using the \texttt{IOHanalyzer} tool \cite{IOHanalyzer}. The heatmaps are generated according to the fraction of times the strategies' mean values are better over all the considered functions in 64D -- with red reflecting total dominance versus blue reflecting total underperformance. 
    The calculations are conducted in isolation for both the separable (LEFT) and non-separable (RIGHT) subsets. 
    The strategies are ordered according to their obtained ranking, so the reader must give mind to the associated labels per each row/column. 
    Importantly, the \texttt{uncorrDG} strategy dominates the separable subset with the \texttt{corrDG} as a first runner-up, whereas their ranks are flipped when shifting to the non-separable subset -- consistently across dimensions (32D and 128D are not shown here). }\label{fig:pairwise_ranking_merged}
\end{figure}

The four main IESs variants, excluding the (1+1) variants, were applied to IQP problems of the form: 
\begin{align}\label{eq:problemClass}
\begin{array}{ll}
\displaystyle  \textrm{minimize}_{\vec{x}} & \frac{1}{c} \cdot \left[ \left( \vec{x}-\vec{\xi}_0 \right)^T \mathbf{H}\left( \vec{x}-\vec{\xi}_0 \right) \right] \\ 
\displaystyle \textrm{subject to:} & \displaystyle \vec{x}\in\mathbb{Z}^{n}.
\end{array}
\end{align}
where the Hessian matrix $\mathbf{H}$, its parametric condition number $c$ and the location vector $\vec{\xi}_0$ completely define a problem instance. 
To instantiate a testbed, the commonly used Ellipsoidal, Discus and Cigar functions \cite{BBOB_GECCO2010}\footnote{These functions correspond to $\left\{f_2,~f_{10},~f_{11},~f_{12}\right\}$ at the renowned BBOB suite. The actual rotation of $f_{10}$ is implemented as reported in \cite{Shir-Yehudayoff_TCS2020}.} were considered via four Hessian matrices: 
\begin{enumerate}[nosep]
\item[\textbf{H-1}] \text{Discus}:~$\left( \mathcal{H}_{\textrm{disc}}\right)_{11} = c,~\left( \mathcal{H}_{\textrm{disc}}\right)_{ii} = 1 \quad i=2,\ldots, n$;
\item[\textbf{H-2}] \text{Cigar}:~$\left( \mathcal{H}_{\textrm{cigar}}\right)_{11} = 1,~\left( \mathcal{H}_{\textrm{cigar}}\right)_{ii} = c \quad i=2,\ldots,n$;
\item[\textbf{H-3}] Rotated Ellipse: $\mathcal{H}_{\textrm{RE}} = \mathcal{O}\mathcal{H}_{\textrm{ellipse}} \mathcal{O}^{-1}$, where $\mathcal{O}$ is rotation by $\approx \frac{\pi}{4}$ radians in the plane spanned by $(1,0,1,0,\ldots)^T$ and $(0,1,0,1,\ldots)^T$;
\item[\textbf{H-4}] Hadamard Ellipse: $\mathcal{H}_{\textrm{HE}} = \mathcal{S}\mathcal{H}_{\textrm{ellipse}} \mathcal{S}^{-1}$, where the rotation constitutes the normalized Hadamard matrix, $\mathcal{S}:=\textrm{Hadamard}(n)/\sqrt{n}$.
\end{enumerate}
providing two separable and two non-separable scalable test-cases.
Six levels of conditioning, $c\in \left\{10,10^2,\ldots,10^6 \right\}$ yielded altogether 24 problem instances per dimensionality.
Systematic benchmarking of the IESs over this testbed included three search-space dimensionalities $n \in \{32,64,128 \}$. 
 Figure \ref{fig:pairwise_ranking_merged} presents pairwise comparisons of the four IESs via heatmaps, according to the fraction of times their mean values are better over all the functions considered, per $n=64$, where the calculations were conducted in isolation for both the separable (left) and non-separable subsets (right).
This pattern of behavior generalized to 32D and 128D, and enabled that study to reach sharp comparative performance conclusions with respect to this testbed: 
\begin{itemize}[nosep]
    \item The \text{uncorrelatedDG} IES dominates the separable subset; the \texttt{correlatedDG} IES is always the first runner-up.
    \item The \texttt{correlatedDG} IES dominates the nonseparable subset; the \texttt{uncorrelatedDG} IES is always the first runner-up.
    \item  DG-based IESs always outperform TN-based IESs over the tested integer suite.
\end{itemize}

\section{Summary and Conclusions} \label{sec:summary}
This study examined a dedicated mutation operator for integer variables in evolution strategies.
We first clarified both the conceptual distance—and the common ground—between integer and continuous unconstrained optimization from an evolutionary-algorithms perspective. 
We then pointed out that even in quadratic unconstrained minimization with a positive-definite form matrix (which yields a convex optimization problem in the continuous case), multimodality can arise when there are dependencies between the variables, i.e., non-zero entries on the diagonal of the form matrix. 
Central to this discussion is our proposed shift of attention to $\ell_{1}$-norm symmetries: whereas classical analyses typically assume $\ell_{2}$-invariance, we showed that the combinatorial setting demands symmetry with respect to Manhattan distances. 
Even a convex quadratic objective can exhibit multimodality once its domain is discretized to the integer lattice, underscoring the need for mutation operators that respect the geometry of the discrete space.

Building on this insight, we proposed a procedure for drawing correlated integer samples that exactly follow the DG distribution, as proposed for independent sampling by Rudolph \cite{rudolph1994evolutionary}. 
Visual inspections and correlation statistics confirm that the intended dependence structure is faithfully reproduced. 
At the same time, it is currently unclear how effectively the considered statistical measures can express an $\ell_1$-norm-driven correlation in a meaningful sense -- it requires more investigation and is reserved for future research. 
Moreover, a comparative study of candidate mutation kernels—including the TN and the DG distributions -- highlighted the latter’s favorable entropy profile; among all tested distributions it delivers also in the correlated case the largest entropy for a given expected step length, thus offering the richest exploratory power.

Empirical tests on discretized benchmark functions reveal significantly better convergence when the IES employs the DG-based mutation operator. 
Yet every configuration, irrespective of the mutation kernel, experiences an apparently unavoidable loss of accuracy once the search enters the immediate vicinity of the optimum.  
The resulting stagnation sets in almost deterministically at a characteristic distance, reminiscent of a \textit{phase transition}, while simple tweaks fail to remove it.

Looking ahead, several avenues merit closer attention.  
A rigorous run-time analysis, especially of the final search phase, could uncover the mechanisms behind the observed stagnation and inspire strategies that maintain progress to exact optima.  
Because the DG law generates step lengths in closed form, it is a natural candidate for derandomized step-size adaptation schemes such as the CMA-ES; combining the DG distribution might even further improve results.  
Further experiments on problems with correlated integer decision variables -- NK landscapes \cite{li2013mixed} provide an immediate MI testbed, as well as real-world instances -- will demonstrate how well the method generalizes to MI cases. In the MI case it might however be difficult to capture the dependencies amone integer and continuous variables \cite{emmerich2008mixed}.
Finally, enriching the algorithm with boundary-aware mutations should enable effective search when true optima lie on, or close to, integer constraints, and extending the framework to MI settings will widen its practical reach.

\end{document}